\documentclass[12pt]{article}
\usepackage{amsmath,amssymb}
\usepackage{theorem}
\usepackage{mathdots}
\usepackage{graphicx}
\usepackage{xypic}
\usepackage{epic,eepic}
\usepackage{hyperref}
\usepackage{doi}

\newtheorem{theorem}{Theorem}[section]
\newtheorem{proposition}[theorem]{Proposition}
\newtheorem{lemma}[theorem]{Lemma}
\newtheorem{pfpf}{{\it Proof.}}

\newenvironment{prf}{\begin{pfpf}\rm}{\hspace*{\fill}{$\square$}\end{pfpf}}
\theorembodyfont{\rmfamily}\newtheorem{remark}[theorem]{Remark}
\theorembodyfont{\rmfamily}\newtheorem{example}{Example}
\theorembodyfont{\rmfamily}\newtheorem{definition}[theorem]{Definition}

\newcommand{\InjHull}{\operatorname{E}}
\newcommand{\tr}{\operatorname{tr}}

\numberwithin{equation}{section}

\begin{document}

\title{\textbf{Injective hulls of odd cycles}\\[9mm]}

\author{\Large Ruedi Suter\\[3mm]
{\textit{Department of Mathematics}}\\
{\textit{ETH Zurich}}\\
{\textit{Raemistrasse 101, 8092 Zurich, Switzerland}}\\[2mm]
e-mail: \texttt{suter@math.ethz.ch}}

\maketitle

\begin{abstract}
The injective hulls of odd cycles are described explicitly.
\end{abstract}

\section{Introduction}
Let $\lambda=(\lambda_1,\lambda_2,\dots,\lambda_m)$ and
$\mu=(\mu_1,\mu_2,\dots,\mu_n)$ be integer partitions, that is, weakly
decreasing finite sequences of positive integers. Their distance
$d(\lambda,\mu)$ is the distance that is induced from the Hasse diagram
of Young's lattice when all edges have length one. In other words,
$d(\lambda,\mu)$ is the number of boxes in the Young diagrams of
$\lambda$ and $\mu$ that belong to $\lambda$ but not to $\mu$ or vice
versa. Hence
$$d(\lambda,\mu)=\sum_{l=1}^m\lambda_l+\sum_{l=1}^n\mu_l-
2\sum_{l=1}^{\min\{m,n\}}\min\bigl\{\lambda_l,\mu_l\bigr\}=|\lambda|+|\mu|
-2|\lambda\cap\mu|$$
where absolute signs denote the size of a partition,
and $\lambda\cap\mu$ is the partition whose Young diagram is, with an evident
meaning, the  intersection of the Young diagrams of $\lambda$ and $\mu$.

For $N\geqslant1$ let $\mathbb Y_N$ be the set of those integer
partitions that have maximal hook length ($\lambda_1+m-1$ for $\lambda\neq()$
as above; the empty partition has maximal hook length zero) strictly
smaller than $N$. 
The $N$ rectangular partitions $R_j=(j^{N-j})$ for $j=0,\dots,N-1$,
where $(0^N)$ is synonymous with the empty partition $()$,
belong to $\mathbb Y_N$. The $N-1$ partitions $R_1,\dots,R_{N-1}$ are the
maximal elements in $\mathbb Y_N$ (for $N\geqslant2$), and the empty partition
$R_0$ is of course the least element. The cardinality of $\mathbb Y_N$ is
$2^{N-1}$. 

We consider the $N$-point metric space $X_N=\{R_0,\dots,R_{N-1}\}$.
The distances are
\begin{align}\label{distR}
d(R_i,R_j)=|j-i|\bigl(N-|j-i|\bigr).
\end{align}
It is convenient to have $R_N=()$ as an alternative notation for the empty
partition $R_0$, so that the formula (\ref{distR}) for the distances between
$R_i$ and $R_j$ holds for $i,j\in\{0,\dots,N\}$. The symmetric matrix
$\bigl(d(R_i,R_j)\bigr)_{i,j=0,\dots,N-1}$ is circulant. In fact, much more is
true. There is a cyclic action of order $N$ on $\mathbb Y_N$ by isometries.
The cyclic action $\tau:\mathbb Y_N\to\mathbb Y_N$ is given by the formula
$$\lambda=(\lambda_1,\dots,\lambda_m)\stackrel{\tau}{\longmapsto}
\tau(\lambda)=(N-m-1,\lambda_1-1,\dots,\lambda_m-1)|_{\textup{remove trailing zeros}}$$
(for the empty partition put $m=0$, $\lambda_1=1$). In particular, we get
$\tau(R_{j+1})=R_j$ for $j=0,\dots,N-1$. It is straightforward to check
that $\tau$ preserves distances. For $\lambda$ and $\mu$ as above with
$n\leqslant m$ we compute
\begin{align*}
d\bigl(\tau(\lambda),\tau(\mu)\bigr)&=\bigl|\tau(\lambda)\bigr|
+\bigl|\tau(\mu)\bigr|-2\bigl|\tau(\lambda)\cap\tau(\mu)\bigr|\\[3mm]
&=N-m-1+|\lambda|-m+N-n-1+|\mu|-n\\[-1mm]
&\phantom{=}\quad{}-2\Bigl(N-m-1+\sum_{l=1}^n
\min\bigl\{\lambda_l-1,\mu_l-1\bigr\}\Bigr)\\
&=|\lambda|+|\mu|-2\sum_{l=1}^n\min\bigl\{\lambda_l,\mu_l\bigr\}
=d(\lambda,\mu).
\end{align*}

The injective hull $\InjHull(X_N)$ of the metric space $X_N$ was described
in \cite{Su2}. The $1$-skeleton of $\InjHull(X_N)$ can be identified with
the Hasse graph of the subposet of Young's lattice that is induced from the
vertex set $\mathbb Y_N$. In particular, this yielded a new geometric
interpretation of the cyclic symmetries described above. The original
combinatorial approach is in \cite{Su1}, and another geometric description
is mentioned in \cite{Su1a}.

$\InjHull(X_N)$ is an example of a most complicated injective hull of an
$N$-point metric space in the sense that it realizes the maximal possible
number of $v$-faces that such an injective hull can have. On the other hand,
it is very special in that it exhibits an $N$-fold cyclic symmetry
(more precisely, the symmetry group is a dihedral group of order $2N$,
at least if $N\geqslant3$).

\section{$N$-cycle metric spaces}
For $N\geqslant4$ the $N$-point metric space $X_N=\{R_0,\dots,R_{N-1}\}$
is not a cycle because $d(R_0,R_2)=2(N-2)<2(N-1)=d(R_0,R_1)+d(R_1,R_2)$.
But an $N$-cycle metric space can be realized as a subspace of $\mathbb Y_N$.
For this let $k=\bigl\lfloor\frac N2\bigr\rfloor$ and consider the orbit of
the staircase partition $\alpha_0:=(k-1,\stackrel{\searrow}{\ldots},1)
\in\mathbb Y_N$, namely
$C_N:=\bigl\{\tau^j(\alpha_0)\bigm|j=0,\dots,N-1\bigr\}$.

For $N=2k$ even
$$\alpha_j:=\tau^j(\alpha_0)=\begin{cases}
(k-1,\stackrel{\searrow}{\ldots},1)&\mbox{for $j=0$ (and $j=N$)}\\
(k,\stackrel{\searrow}{\ldots},k-j+1,k-j-1,\stackrel{\searrow}{\ldots},1)
&\mbox{for $j=1,\dots,k-2$}\\
(k,\stackrel{\searrow}{\ldots},2)&\mbox{for $j=k-1$}\\
(k,\stackrel{\searrow}{\ldots},1)&\mbox{for $j=k$}\\
(k-1,\stackrel{\searrow}{\ldots},N-j,N-j,\stackrel{\searrow}{\ldots},1)
&\mbox{for $j=k+1,\dots,N-1$}
\end{cases}$$
where the partition $\alpha_j$ has $k-1$ parts for $j=0,\dots,k-1$ and
$k$ parts for $j=k,\dots,N-1$.
Notice that $d(\alpha_i,\alpha_j)=\min\{|j-i|,N-|j-i|\}$
(it suffices to verify that $d(\alpha_0,\alpha_1)=1$ and
$d(\alpha_0,\alpha_k)=k$).

For $N=2k+1$ odd
$$\alpha_j:=\tau^j(\alpha_0)=\begin{cases}
(k-1,\stackrel{\searrow}{\ldots},1)&\mbox{for $j=0$ (and $j=N$)}\\
(k+1,\stackrel{\searrow}{\ldots},k-j+2,k-j-1,\stackrel{\searrow}{\ldots},1)&\mbox{for $j=1,\dots,k-2$}\\
(k+1,\stackrel{\searrow}{\ldots},3)&\mbox{for $j=k-1$}\\
(k+1,\stackrel{\searrow}{\ldots},2)&\mbox{for $j=k$}\\
(k,k,\stackrel{\searrow}{\ldots},1)&\mbox{for $j=k+1$}\\
(k-1,\stackrel{\searrow}{\ldots},N-j,N-j,N-j,\stackrel{\searrow}{\ldots},1)&\mbox{for $j=k+2,\dots,N-1$}
\end{cases}$$
where the partition $\alpha_j$ has $k-1$ parts for $j=0,\dots,k-1$ and
$\alpha_k$ has $k$ parts and $\alpha_j$ has $k+1$ parts for $j=k+1,\dots,N-1$.
Observe that $d(\alpha_i,\alpha_j)=2\min\{|j-i|,N-|j-i|\}$
(it suffices to verify that $d(\alpha_0,\alpha_1)=2$ (if $k>0$) and
$d(\alpha_0,\alpha_k)=2k$).

In other words, for $N$ of either parity
\begin{align}\label{metricevenodd}
d(\alpha_i,\alpha_j)=\min\bigl\{|j-i|,N-|j-i|\bigr\}\cdot d(\alpha_0,\alpha_1)
\end{align}
and the formula is also valid when $i=N$ or $j=N$ (with $\alpha_N=\alpha_0$). 
This means that the metric space $C_N$ is an $N$-cycle.
Since $C_N$ embeds isometrically into $\InjHull(X_N)$, the smallest
injective subspace of $\InjHull(X_N)$ that contains the image of $C_N$
can be taken as a model for the injective hull of $C_N$.

For $N$ even the injective hull of $C_N$ can be seen as the central
$\frac N2$-dimensional cube that was mentioned in the conclusions-and-outlook
section in \cite{Su2}. The fact that the injective hull of a $2k$-cycle is
just a $k$-cube is the content of Section~9 in \cite{GM}.
For $N$ odd the situation is somewhat more involved.

\section{Discrete M\"obius strips and integer partitions}
We need to recall some notations from \cite{Su2}. The discrete
M\"obius strip $\mathfrak X_N$ has $\frac12 N(N+1)$ sites
$(i,j)$ for $0\leqslant i\leqslant j\leqslant N$, where
$(i,N)=(0,i)\in\mathfrak X_N$ for $i=0,\dots,N$.
The following picture (for $N=9$) illustrates the setup. We have
$O=(0,0)$, $O'=(0,N)$, $O''=(N,N)$ and $O=O'=O''$ as sites in $\mathfrak X_N$;
$P=(0,2)$, $P'=(2,N)$; $Q=(0,3)$, $Q'=(3,N)$; $R=(0,N-1)$, $R'=(N-1,N)$;
and $P=P'$, $Q=Q'$, $R=R'$ as sites in the M\"obius strip. For each partition
$\lambda\in\mathbb Y_N$ we get the corresponding outer rim $\mathcal L_\lambda
\subseteq\mathfrak X_N$, which is a homotopically nontrivial loop consisting
of $N$ sites in the M\"obius strip. In the picture we have
$\lambda=(5,5,2,1)=(5^2\,2\,1)$, whose Young diagram, displayed in the Russian
convention, appears
below the outer rim $\mathcal L_\lambda$.
\enlargethispage*{5mm}
\setlength{\unitlength}{0.0005in}
\begin{center}
\begin{picture}(6024,3339)(0,-10)
\path(12,3012)(3012,12)(6012,3012)
	(5712,3312)(2712,312)
\path(5712,2712)(5112,3312)(2412,612)
\path(5112,2112)(3912,3312)(1812,1212)
\path(4812,1812)(3312,3312)(1512,1512)
\path(4512,1512)(2712,3312)(1212,1812)
\path(4212,1212)(2112,3312)(912,2112)
\path(3312,312)(312,3312)(12,3012)
\path(3912,912)(1512,3312)(612,2412)
\path(3612,612)(912,3312)(312,2712)
\path(5412,2412)(4512,3312)(2112,912)
\texture{55888888 88555555 5522a222 a2555555 55888888 88555555 552a2a2a 2a555555 
	55888888 88555555 55a222a2 22555555 55888888 88555555 552a2a2a 2a555555 
	55888888 88555555 5522a222 a2555555 55888888 88555555 552a2a2a 2a555555 
	55888888 88555555 55a222a2 22555555 55888888 88555555 552a2a2a 2a555555 }
\shade\path(1512,1512)(2112,2112)(2412,1812)
	(2712,2112)(3012,1812)(3912,2712)
	(4812,1812)(4512,1512)(3912,2112)
	(3012,1212)(2712,1512)(2412,1212)
	(2112,1512)(1812,1212)(1512,1512)
\path(1512,1512)(2112,2112)(2412,1812)
	(2712,2112)(3012,1812)(3912,2712)
	(4812,1812)(4512,1512)(3912,2112)
	(3012,1212)(2712,1512)(2412,1212)
	(2112,1512)(1812,1212)(1512,1512)
\put(312,3012){\makebox(0,0){\footnotesize$O$}}
\put(3012,312){\makebox(0,0){\footnotesize$O'$}}
\put(5712,3012){\makebox(0,0){\footnotesize$O''$}}
\put(912,2412){\makebox(0,0){\footnotesize$P$}}
\put(3612,912){\makebox(0,0){\footnotesize$P'$}}
\put(1212,2112){\makebox(0,0){\footnotesize$Q$}}
\put(3912,1212){\makebox(0,0){\footnotesize$Q'$}}
\put(2712,612){\makebox(0,0){\footnotesize$R$}}
\put(5412,2712){\makebox(0,0){\footnotesize$R'$}}
\put(4900,1540){\makebox(0,0){$\mathcal L_\lambda$}}
\end{picture}
\end{center}
The empty partition has as its outer rim
$\mathcal L_{()}=\bigl\{(0,0),\dots,(0,N),\dots,(N,N)\bigr\}$,
where all its $2N+1$ positions in the triangular shape are listed, but
$\mathcal L_{()}$ has only $N$ sites in the M\"obius strip.
The mapping $\lambda\mapsto\mathcal L_\lambda$ is a bijection between
$\mathbb Y_N$ and the set of homotopically nontrivial loops of length
$N$ in the M\"obius strip $\mathfrak X_N$. The transported action $\tau$
comes from translating the loops in the M\"obius strip, for instance
in the example above observe the outer rims of the partitions
$(6,3,2)\stackrel{\tau}{\longmapsto}(5,5,2,1)\stackrel{\tau}{\longmapsto}
(4,4,4,1)$.

For $N=2k+1$ odd we consider the discrete M\"obius strip
$\mathfrak X_N^\circ\subseteq\mathfrak X_N$ that is defined by
$$\mathfrak X_N^\circ:=\bigl\{(i,j)\in\mathfrak X_N\bigm|
d(\alpha_i,\alpha_j)\geqslant2(k-1)\bigr\}.$$
Its $2N$ sites are
\begin{align*}
&(0,k-1),(1,k),\dots,(k+1,N-1),(0,k+2),(1,k+3),\dots,(k-2,N-1),\\
&(0,k),(1,k+1),\dots,(k,N-1),(0,k+1),(1,k+2),\dots,(k-1,N-1).
\end{align*}

\begin{example}
For $N=9$ the sites in $\mathfrak X_N$ that are
crossed \mbox{\LARGE$\times$} in the picture lie outside
$\mathfrak X_N^\circ$. The positions in the triangular shape are
$$\begin{array}{llll}
A=(0,k-1)&B=(0,k)&C=(0,k+1)&D=(0,k+2)\\
A'=(k-1,N)&B'=(k,N)&C'=(k+1,N)&D'=(k+2,N)
\end{array}$$
and $A=A'$, $B=B'$, $C=C'$, and $D=D'$ as sites in $\mathfrak X_N^\circ$.
\setlength{\unitlength}{0.0005in}
\begin{center}
\begin{picture}(6024,3339)(0,-10)
\path(12,3012)(3012,12)(6012,3012)
	(5712,3312)(2712,312)
\path(5712,2712)(5112,3312)(2412,612)
\path(5412,2412)(4512,3312)(2112,912)
\path(5112,2112)(3912,3312)(1812,1212)
\path(4812,1812)(3312,3312)(1512,1512)
\path(4512,1512)(2712,3312)(1212,1812)
\path(4212,1212)(2112,3312)(912,2112)
\path(3912,912)(1512,3312)(612,2412)
\path(3612,612)(912,3312)(312,2712)
\path(3312,312)(312,3312)(12,3012)
\put(312,3012){\makebox(0,0){\LARGE$\times$}}
\put(912,3012){\makebox(0,0){\LARGE$\times$}}
\put(1512,3012){\makebox(0,0){\LARGE$\times$}}
\put(2112,3012){\makebox(0,0){\LARGE$\times$}}
\put(2712,3012){\makebox(0,0){\LARGE$\times$}}
\put(3312,3012){\makebox(0,0){\LARGE$\times$}}
\put(3912,3012){\makebox(0,0){\LARGE$\times$}}
\put(4512,3012){\makebox(0,0){\LARGE$\times$}}
\put(5112,3012){\makebox(0,0){\LARGE$\times$}}
\put(5712,3012){\makebox(0,0){\LARGE$\times$}}
\put(612,2712){\makebox(0,0){\LARGE$\times$}}
\put(1212,2712){\makebox(0,0){\LARGE$\times$}}
\put(1812,2712){\makebox(0,0){\LARGE$\times$}}
\put(2412,2712){\makebox(0,0){\LARGE$\times$}}
\put(3012,2712){\makebox(0,0){\LARGE$\times$}}
\put(3612,2712){\makebox(0,0){\LARGE$\times$}}
\put(4212,2712){\makebox(0,0){\LARGE$\times$}}
\put(4812,2712){\makebox(0,0){\LARGE$\times$}}
\put(5412,2712){\makebox(0,0){\LARGE$\times$}}
\put(912,2412){\makebox(0,0){\LARGE$\times$}}
\put(1512,2412){\makebox(0,0){\LARGE$\times$}}
\put(2112,2412){\makebox(0,0){\LARGE$\times$}}
\put(2712,2412){\makebox(0,0){\LARGE$\times$}}
\put(3312,2412){\makebox(0,0){\LARGE$\times$}}
\put(3912,2412){\makebox(0,0){\LARGE$\times$}}
\put(4512,2412){\makebox(0,0){\LARGE$\times$}}
\put(5112,2412){\makebox(0,0){\LARGE$\times$}}
\put(2412,912){\makebox(0,0){\LARGE$\times$}}
\put(3012,912){\makebox(0,0){\LARGE$\times$}}
\put(3612,912){\makebox(0,0){\LARGE$\times$}}
\put(2712,612){\makebox(0,0){\LARGE$\times$}}
\put(3312,612){\makebox(0,0){\LARGE$\times$}}
\put(3012,312){\makebox(0,0){\LARGE$\times$}}
\put(1212,2112){\makebox(0,0){\footnotesize$A$}}
\put(1512,1812){\makebox(0,0){\footnotesize$B$}}
\put(1812,1512){\makebox(0,0){\footnotesize$C$}}
\put(2112,1212){\makebox(0,0){\footnotesize$D$}}
\put(3912,1212){\makebox(0,0){\footnotesize$A'$}}
\put(4212,1512){\makebox(0,0){\footnotesize$B'$}}
\put(4512,1812){\makebox(0,0){\footnotesize$C'$}}
\put(4812,2112){\makebox(0,0){\footnotesize$D'$}}
\end{picture}
\end{center}
\end{example}

\enlargethispage*{8mm}
\begin{definition}\label{YNcirc}
Let
$$\mathbb Y_N^\circ:=
\bigl\{\lambda\in\mathbb Y_N\bigm|\mathcal L_\lambda\subseteq
\mathfrak X_N^\circ\bigr\}$$
be the set of those partitions with maximal hook length ${}<N$ whose outer
rim is contained in $\mathfrak X_N^\circ$. Note that $\alpha_0\in
\mathbb Y_N^\circ$. The cyclic action $\tau$ restricts to
$\tau:\mathbb Y_N^\circ\to\mathbb Y_N^\circ$.
\end{definition}

For $(i,j)\in\mathfrak X_N^\circ$ we have according to (\ref{distR}) and
(\ref{metricevenodd})
$$d(R_i,R_j)=\begin{cases}
k^2+k&\mbox{if $d(\alpha_i,\alpha_j)=2k$,}\\
k^2+k-2&\mbox{if $d(\alpha_i,\alpha_j)=2k-2$,}
\end{cases}$$
or in other words
\begin{align}\label{verschieb}
d(R_i,R_j)=d(\alpha_i,\alpha_j)+k^2-k\qquad
\mbox{for $(i,j)\in\mathfrak X_N^\circ$}.
\end{align}

\section{Injective hulls}
The injective hulls of the $N$-point metric spaces $X_N$ and $C_N$ can
be realized as the  polyhedral complexes $\InjHull(X_N)$ and
$\InjHull(C_N)$ defined as follows:
\begin{align}\notag
\Delta(X_N)&=\bigl\{f\in\mathbb R^{X_N}\bigm|
\forall R_i,R_j\in X_N: f(R_i)+f(R_j)\geqslant d(R_i,R_j)\bigr\}\\\label{hullX}
\InjHull(X_N)&=\bigl\{f\in\Delta(X_N)\bigm|
\forall R_i\in X_N\,\exists R_j\in X_N:f(R_i)+f(R_j)=d(R_i,R_j)\bigr\}
\intertext{and}\notag
\Delta(C_N)&=\bigl\{g\in\mathbb R^{C_N}\bigm|
\forall\alpha_i,\alpha_j\in C_N: g(\alpha_i)+g(\alpha_j)\geqslant d
(\alpha_i,\alpha_j)\bigr\}\\\label{hullC}
\InjHull(C_N)&=\bigl\{g\in\Delta(C_N)\bigm|
\forall\alpha_i\in C_N\,\exists\alpha_j\in C_N:g(\alpha_i)+g(\alpha_j)
=d(\alpha_i,\alpha_j)\bigr\}.
\end{align}
We use the bijection $\iota:C_N\to X_N$, $\alpha_j\mapsto R_j$ to map
$\mathbb R^{X_N}\ni f\mapsto f\circ\iota\in\mathbb R^{C_N}$.

The injective hull $\InjHull(X_N)$ was studied in \cite{Su2}. Its
vertices are parametrized by the partitions in $\mathbb Y_N$:
\begin{align}\notag
\lambda\in\mathbb Y_N&\rightsquigarrow\mathcal L_\lambda\subseteq
\mathfrak X_N&&\mbox{its outer rim}\\\notag
&\rightsquigarrow f_\lambda\in\InjHull(X_N)
&&\mbox{the solution of $f(R_{i_l})+f(R_{j_l})=d(R_{i_l},R_{j_l})$
($l=0,\dots,N-1$)}\\\label{systemregular}
&&&\mbox{where $\mathcal L_\lambda=\bigl\{(i_l,j_l)\bigm|l=0,\dots,N-1\bigr\}$,}\\\notag
&&&\mbox{explicitly: $f_\lambda(R_j)=\bigl|\tau^j(\lambda)\bigr|$}.
\end{align}

If $\lambda\in\mathbb Y_N$ has $s$ inner corners and $\lambda_\downarrow$
denotes the partition that is got from $\lambda$ by removing all its
inner corners, then the convex hull of the $2^s$ vertices $f_\nu$ where
$\nu\in[\lambda_\downarrow,\lambda]$ is a subset of $\InjHull(X_N)$.
Thus there are $\binom sv$ $v$-faces with ``top vertex'' $f_\lambda$.

Let again $N=2k+1$. The constant function $o\in\mathbb R^{C_N}$
is defined as
$o(\alpha_j)=\frac12(k^2-k)=f_{\alpha_0}(R_0)$ ($j=0,\dots,N-1$).
Using (\ref{verschieb}), we obtain for $\lambda\in\mathbb Y_N^\circ$
by taking $f_\lambda\in\InjHull(X_N)$ and comparing (\ref{hullX}) with
(\ref{hullC})  
\begin{align}\label{glambda}
g_\lambda:=f_\lambda\circ\iota-o\in\InjHull(C_N).
\end{align}
In the same way, for $f\in\InjHull(X_N)$ in a $v$-face with vertices
$f_\nu$ with all $\nu\in\mathbb Y_N^\circ$, we have 
$g:=f\circ\iota-o\in\InjHull(C_N)$.
 
Let $\InjHull(X_N)^\circ$ be the
subcomplex of $\InjHull(X_N)$ induced from the vertex set $f_\lambda$ for
$\lambda\in\mathbb Y_N^\circ$. Then $\InjHull(X_N)^\circ$  can be identified
via the translation $f\mapsto f\circ\iota-o$ with a subset of $\InjHull(C_N)$.
It turns out that the image of $\InjHull(X_N)^\circ$ is actually the whole
of $\InjHull(C_N)$. To prove this, we use a $1$-Lipschitz retraction
$\InjHull(X_N)\to\InjHull(X_N)^\circ$.

\setlength{\unitlength}{0.0006in}
\begin{center}
\begin{picture}(7224,5549)(0,-490)
\path(5862,3622)(4362,5122)(2112,2872)
\path(5712,3472)(4062,5122)(1962,3022)
\path(6612,4372)(5862,5122)(2862,2122)
\path(5562,3322)(3762,5122)(1812,3172)
\path(7062,4822)(6762,5122)(3312,1672)
\path(6762,4522)(6162,5122)(3012,1972)
\path(6462,4222)(5562,5122)(2712,2272)
\path(6162,3922)(4962,5122)(2412,2572)
\path(6012,3772)(4662,5122)(2262,2722)
\path(5412,3172)(3462,5122)(1662,3322)
\path(5262,3022)(3162,5122)(1512,3472)
\path(5112,2872)(2862,5122)(1362,3622)
\path(4962,2722)(2562,5122)(1212,3772)
\path(4812,2572)(2262,5122)(1062,3922)
\path(4662,2422)(1962,5122)(912,4072)
\path(4062,1822)(762,5122)(312,4672)
\path(3912,1672)(462,5122)(162,4822)
\path(3762,1522)(162,5122)(12,4972)
\path(6312,4072)(5262,5122)(2562,2422)
\path(4362,2122)(1362,5122)(612,4372)
\path(12,4972)(3612,1372)(7212,4972)
	(7062,5122)(3462,1522)
\path(6912,4672)(6462,5122)(3162,1822)
\path(4212,1972)(1062,5122)(462,4522)
\path(1662,472)(1812,622)(2262,172)
	(2712,622)(3162,172)(3612,622)
	(4062,172)(4512,622)(4962,172)(5412,622)
\path(1812,322)(2112,622)(2562,172)
	(3012,622)(3462,172)(3912,622)
	(4362,172)(4812,622)(5262,172)
\path(1962,172)(2412,622)(2862,172)
	(3312,622)(3762,172)(4212,622)
	(4662,172)(5112,622)(5412,322)
\path(5412,622)(5562,472)
\path(4512,2272)(1662,5122)(762,4222)
\texture{55888888 88555555 5522a222 a2555555 55888888 88555555 552a2a2a 2a555555 
	55888888 88555555 55a222a2 22555555 55888888 88555555 552a2a2a 2a555555 
	55888888 88555555 5522a222 a2555555 55888888 88555555 552a2a2a 2a555555 
	55888888 88555555 55a222a2 22555555 55888888 88555555 552a2a2a 2a555555 }
\shade\path(1662,3322)(2262,3922)(2412,3772)
	(2862,4222)(3762,3322)(3912,3472)
	(4212,3172)(4362,3022)(4512,3172)
	(4662,3322)(4812,3172)(5112,3472)
	(5412,3172)(5262,3022)(5112,3172)
	(4812,2872)(4662,3022)(4362,2722)
	(3912,3172)(3762,3022)(2862,3922)
	(2412,3472)(2262,3622)(1812,3172)(1662,3322)
\path(1662,3322)(2262,3922)(2412,3772)
	(2862,4222)(3762,3322)(3912,3472)
	(4212,3172)(4362,3022)(4512,3172)
	(4662,3322)(4812,3172)(5112,3472)
	(5412,3172)(5262,3022)(5112,3172)
	(4812,2872)(4662,3022)(4362,2722)
	(3912,3172)(3762,3022)(2862,3922)
	(2412,3472)(2262,3622)(1812,3172)(1662,3322)
\shade\path(1662,472)(1962,772)(2112,622)
	(2262,772)(2412,622)(2562,772)
	(2712,622)(2862,772)(3012,622)
	(3162,772)(3312,622)(3462,772)
	(3762,472)(3912,622)(4212,322)
	(4362,472)(4512,322)(4662,472)
	(4812,322)(5112,622)(5412,322)
	(5262,172)(5112,322)(4812,22)
	(4662,172)(4512,22)(4362,172)
	(4212,22)(3912,322)(3762,172)
	(3462,472)(3312,322)(3162,472)
	(3012,322)(2862,472)(2712,322)
	(2562,472)(2412,322)(2262,472)
	(2112,322)(1962,472)(1812,322)(1662,472)
\path(1662,472)(1962,772)(2112,622)
	(2262,772)(2412,622)(2562,772)
	(2712,622)(2862,772)(3012,622)
	(3162,772)(3312,622)(3462,772)
	(3762,472)(3912,622)(4212,322)
	(4362,472)(4512,322)(4662,472)
	(4812,322)(5112,622)(5412,322)
	(5262,172)(5112,322)(4812,22)
	(4662,172)(4512,22)(4362,172)
	(4212,22)(3912,322)(3762,172)
	(3462,472)(3312,322)(3162,472)
	(3012,322)(2862,472)(2712,322)
	(2562,472)(2412,322)(2262,472)
	(2112,322)(1962,472)(1812,322)(1662,472)
\put(5862,3022){\makebox(0,0)[lb]{$\mathcal L_\lambda\subseteq\mathfrak X_{23}$}}
\put(10,1400){\makebox(0,0)[lb]{$\lambda=(11\,9\,7^3\,6^6\,3)\in\mathbb Y_{23}$}}
\put(5862,172){\makebox(0,0)[lb]{$\mathcal L_{\lambda^\circ}
\subseteq\mathfrak X_{23}^\circ$}}
\put(10,-500){\makebox(0,0)[lb]{$\lambda^\circ
=(11\,9\,8\,7^2\,6^2\,5\,4\,3\,2\,1)
\in\mathbb Y_{23}^\circ$}}
\thicklines
\path(1512,622)(2112,22)(2262,172)
	(2412,22)(2562,172)(2712,22)
	(2862,172)(3012,22)(3162,172)
	(3312,22)(3462,172)(3612,22)
	(3762,172)(3912,22)(4062,172)
	(4212,22)(4362,172)(4512,22)
	(4662,172)(4812,22)(4962,172)
	(5112,22)(5712,622)(5562,772)
	(5412,622)(5262,772)(5112,622)
	(4962,772)(4812,622)(4662,772)
	(4512,622)(4362,772)(4212,622)
	(4062,772)(3912,622)(3762,772)
	(3612,622)(3462,772)(3312,622)
	(3162,772)(3012,622)(2862,772)
	(2712,622)(2562,772)(2412,622)
	(2262,772)(2112,622)(1962,772)
	(1812,622)(1662,772)(1512,622)
\path(1512,3472)(2112,2872)(2262,3022)
	(2412,2872)(2562,3022)(2712,2872)
	(2862,3022)(3012,2872)(3162,3022)
	(3312,2872)(3462,3022)(3612,2872)
	(3762,3022)(3912,2872)(4062,3022)
	(4212,2872)(4362,3022)(4512,2872)
	(4662,3022)(4812,2872)(4962,3022)
	(5112,2872)(5712,3472)(5562,3622)
	(5412,3472)(5262,3622)(5112,3472)
	(4962,3622)(4812,3472)(4662,3622)
	(4512,3472)(4362,3622)(4212,3472)
	(4062,3622)(3912,3472)(3762,3622)
	(3612,3472)(3462,3622)(3312,3472)
	(3162,3622)(3012,3472)(2862,3622)
	(2712,3472)(2562,3622)(2412,3472)
	(2262,3622)(2112,3472)(1962,3622)
	(1812,3472)(1662,3622)(1512,3472)
\end{picture}
\end{center}
Let $\lambda\in\mathbb Y_N$ with outer rim $\mathcal L_\lambda$. We obtain a
well-defined outer rim $\mathcal L_{\lambda^\circ}\subseteq \mathfrak X_N^\circ$
by starting with the loop $\mathcal L_\lambda$ and successively folding
sites where the loop turns and that are not in $\mathfrak X_N^\circ$ closer
towards $\mathfrak X_N^\circ$. The picture for the example above
illustrates the result
(after twelve foldings in the upper part and one folding in the lower part
of the triangular shape).

Step by step the folding procedure works as follows (and by the way,
in the proof of Proposition~\ref{nonexpanding} below we make the order 
of steps even more detailed and fold two loops simultaneously).
We assume that $0\leqslant i<j\leqslant N-1$.
If $(i+1,j)\notin\mathfrak X_N^\circ$ belongs to the upper part
(i.\,e.\ $j-i-1\leqslant k-2$), then the site $(i,j+1)$ is closer to
$\mathfrak X_N^\circ$ than $(i+1,j)$. Similarly, if 
$(i,j+1)\notin\mathfrak X_N^\circ$ belongs to the lower part
(i.\,e.\ $j+1-i\geqslant k+3$), then the site $(i+1,j)$ is closer to
$\mathfrak X_N^\circ$ than $(i,j+1)$. 

\setlength{\unitlength}{0.00072in}
\begin{center}
\begin{picture}(6624,2439)(0,-10)
\path(1812,1812)(612,612)(12,1212)
	(1212,2412)(2412,1212)(1812,612)(612,1812)
\path(4812,612)(6012,1812)(6612,1212)
	(5412,12)(4212,1212)(4812,1812)(6012,612)
\put(3312,1212){\makebox(0,0){\LARGE$\leftrightsquigarrow$}}
\put(612,1212){\makebox(0,0){\small$(i,j)$}}
\put(1212,1812){\makebox(0,0){\small$(i+1,j)$}}
\put(1812,1212){\makebox(0,0){\small$(i\!+\!1,j\!+\!1)$}}
\put(4812,1212){\makebox(0,0){\small$(i,j)$}}
\put(6012,1212){\makebox(0,0){\small$(i\!+\!1,j\!+\!1)$}}
\put(5412,612){\makebox(0,0){\small$(i,j+1)$}}
\put(-388,1212){\makebox(0,0){$\mathcal L_{\varphi}\supseteq{}$}}
\put(7012,1212){\makebox(0,0){${}\subseteq\mathcal L_{\psi}$}}
\end{picture}
\end{center}
Note also that in order to be able to finally fold everything
to $\mathfrak X_N^\circ$, we must keep in mind that
$(0,j)=(j,N)$ as sites in $\mathfrak X_N$ (or alternatively we could
use the universal covering of the M\"obius strip, with the additional
benefit that we could neglect the spurious distinction between the
upper and lower parts; however, let us stick to the representation
of the sites of $\mathfrak X_N$ in the triangular shape with corners
$(0,0)$, $(0,N)$, and $(N,N)$).
Clearly, this gives a well-defined idempotent mapping
$\mathbb Y_N\ni\lambda\mapsto\lambda^\circ\in\mathbb Y_N^\circ\subseteq
\mathbb Y_N$. Because the cyclic action is given by translation in the
M\"obius strip, it is evident that the idempotent mapping $\lambda\mapsto
\lambda^\circ$ commutes with the cyclic action:
$\bigl(\tau(\lambda)\bigr)^\circ=\tau(\lambda^\circ)$.

\begin{proposition}\label{nonexpanding}
For $\lambda,\mu\in\mathbb Y_N$ we have $d(\lambda^\circ,\mu^\circ)\leqslant
d(\lambda,\mu)$.
\end{proposition}
\begin{prf}
We construct two sequences
$\lambda=\lambda^{(0)},\dots,\lambda^{(l)}=\lambda^\circ$ and
$\mu=\mu^{(0)},\dots,\mu^{(l)}=\mu^\circ$ that satisfy
$d\bigl(\lambda^{(j+1)},\mu^{(j+1)}\bigr)\leqslant
d\bigl(\lambda^{(j)},\mu^{(j)}\bigr)$ for $j=0,\dots,l-1$. To achieve this,
we scan through $\mathfrak X_N-\mathfrak X_N^\circ$ line by line from
the periphery of $\mathfrak X_N$ inwards, namely going through the site
positions $(1,1),\dots,(N-1,N-1)$ in the upper part, then $(0,N)$ in the
lower part, next $(1,2),\dots,(N-2,N-1)$ in the upper part, then
$(0,N-1),(1,N)$ in the lower part, and so on, finally
$(1,k-1),\dots,(k+2,N-1)$ in the upper part and $(0,k+3),\dots,(k-2,N)$
in the lower part of the triangular shape (again, as an alternative, one could
look at the universal covering).
Suppose our scanning stays at site position $L$ and we are at step $j$.
If $L\in\mathcal L_{\lambda^{(j)}}\cup\mathcal L_{\mu^{(j)}}$, then
move to step $j+1$ as follows:
\begin{itemize}
\item[]
\begin{itemize}
\item[$\circ$]
if $L\in\mathcal L_{\lambda^{(j)}}$ then put 
$\mathcal L_{\lambda^{(j+1)}}:=\bigl(\mathcal L_{\lambda^{(j)}}-\{L\}\bigr)
\cup\{L'\}$ where $L'=L\mp(1,-1)$ with ``$-$'' for $L$ in the upper part
(so that $\lambda^{(j+1)}$ is got from $\lambda^{(j)}$ by removing the inner
corner at $L'$),
and ``$+$'' for $L$ in the lower part (so that $\lambda^{(j+1)}$ is got from
$\lambda^{(j)}$ by adding the outer corner at $L$);
otherwise (that is, if $L\notin\mathcal L_{\lambda^{(j)}}$) let
$\lambda^{(j+1)}:=\lambda^{(j)}$;
\item[$\circ$]
if $L\in\mathcal L_{\mu^{(j)}}$ then put 
$\mathcal L_{\mu^{(j+1)}}:=\bigl(\mathcal L_{\mu^{(j)}}-\{L\}\bigr)
\cup\{L'\}$ where $L'=L\mp(1,-1)$ with ``$-$'' for $L$ in the upper part
(so that $\mu^{(j+1)}$ is got from $\mu^{(j)}$ by removing the inner
corner at $L'$),
and ``$+$'' for $L$ in the lower part (so that $\mu^{(j+1)}$ is got from
$\mu^{(j)}$ by adding the outer corner at $L$);
otherwise (that is, if $L\notin\mathcal L_{\mu^{(j)}}$) let
$\mu^{(j+1)}:=\mu^{(j)}$.
\end{itemize}
If $L'$ was an inner corner (respectively $L$ was an outer corner)
for both $\lambda^{(j)}$ and $\mu^{(j)}$, then
$d\bigl(\lambda^{(j+1)},\mu^{(j+1)}\bigr)=d\bigl(\lambda^{(j)},\mu^{(j)}\bigr)$.
Otherwise,
$d\bigl(\lambda^{(j+1)},\mu^{(j+1)}\bigr)=d\bigl(\lambda^{(j)},\mu^{(j)}\bigr)-1$.
\end{itemize}
If $L\notin\mathcal L_{\lambda^{(j)}}\cup\mathcal L_{\mu^{(j)}}$, scan the next
site till there is no next site to scan, in which case the procedure ends
with $l=j$.
\end{prf}

We get a $1$-Lipschitz retraction $r:\InjHull(X_N)\to\InjHull(X_N)^\circ$
that is defined on the vertices by $r(f_\lambda)=f_{\lambda^\circ}$ and extended
in the obvious way to convex combinations corresponding to higher-dimensional
faces. Namely, such a $v$-face is the convex hull of the $2^v$ vertices
$f_\nu$ for $\nu\in[\lambda_V,\lambda]$ where $\lambda\in\mathbb Y_N$,
$V$ with $|V|=v$ is a subset of the set of inner corners of $\lambda$,
and $\lambda_V$ is the partition that is got from $\lambda$ after removing
the inner corners in $V$. Let $V^\circ\subseteq V$ be the subset of those
inner corners that are still inner corners of $\lambda^\circ$ and whose
removal from $\lambda^\circ$ still gives a partition $\lambda^\circ_{V^\circ}
\in\mathbb Y_N^\circ$ ($V^\circ=V\cap\{\mbox{$\circ$-inner corners of
$\lambda^\circ$}\}$ in the sense of Definition~\ref{circinner}).
The image is the $|V^\circ|$-face which is the convex hull of the
$2^{|V^\circ|}$ vertices $f_\nu$ for $\nu\in[\lambda^\circ_{V^\circ},\lambda^\circ]$.
Loosely speaking, the original $v$-cube in $\InjHull(X_N)$ is collapsed along
the $(V-V^\circ)$-directions, and the resulting $|V^\circ|$-cube is then
moved to $\InjHull(X_N)^\circ$.

Let $h:A\to\InjHull(X_N)^\circ$ be any $1$-Lipschitz mapping from a subspace
$A$ of any metric space $Y$. Composing with the inclusion $\InjHull(X_N)^\circ
\hookrightarrow\InjHull(X_N)$ gives us a $1$-Lipschitz mapping
$\operatorname{incl}\circ h:A\to\InjHull(X_N)$. It extends to a $1$-Lipschitz
mapping $\overline{\operatorname{incl}\circ h}:Y\to\InjHull(X_N)$ because
$\InjHull(X_N)$ is an injective metric space. Finally,
$r\circ\overline{\operatorname{incl}\circ h}:Y\to\InjHull(X_N)^\circ$
extends $h$ as a $1$-Lipschitz mapping. In other words, $\InjHull(X_N)^\circ$
is an injective metric space. Thus
$$\InjHull(C_N)=\bigl\{f\circ\iota-o\bigm|f\in\InjHull(X_N)^\circ\bigr\}.$$

One may wonder what one gets from the following process:
\begin{align*}
\lambda\in\mathbb Y_N&\rightsquigarrow\mathcal L_\lambda\subseteq
\mathfrak X_N&&\mbox{its outer rim}\\
&\rightsquigarrow g_\lambda\in\mathbb R^{C_N}
&&\mbox{the solution of $g(\alpha_{i_l})+g(\alpha_{j_l})
=d(\alpha_{i_l},\alpha_{j_l})$
($l=0,\dots,N-1$)}\\
&&&\mbox{where $\mathcal L_\lambda=\bigl\{(i_l,j_l)\bigm|l=0,\dots,N-1\bigr\}$.}
\end{align*}
Recall that the linear system is regular, as was already implicitly used
in (\ref{systemregular}). If $\lambda\in\mathbb Y_N^\circ$, the functions
$g_\lambda\in\InjHull(C_N)$ were described in (\ref{glambda}).

Suppose the outer rims $\mathcal L_\varphi$ and $\mathcal L_\psi$ are
as before, with $\mathcal L_\varphi\ni(i+1,j)\notin\mathfrak X_N^\circ$
and $j-i-1\leqslant k-2$. Then
\begin{align*}
d(\alpha_i,\alpha_j)=d(\alpha_{i+1},\alpha_{j+1})&=d(\alpha_{i+1},\alpha_j)+2\\
d(\alpha_i,\alpha_{j+1})&=d(\alpha_{i+1},\alpha_j)+4
\end{align*}
so that
$$d(\alpha_i,\alpha_j)+d(\alpha_{i+1},\alpha_{j+1})=
d(\alpha_{i+1},\alpha_j)+d(\alpha_i,\alpha_{j+1})$$
and it follows that under the condition that
\begin{align*}
g(\alpha_i)+g(\alpha_j)&=d(\alpha_i,\alpha_j)\\
g(\alpha_{i+1})+g(\alpha_{j+1})&=d(\alpha_{i+1},\alpha_{j+1})
\end{align*}
the two equations
\begin{align*}
g(\alpha_{i+1})+g(\alpha_j)&=d(\alpha_{i+1},\alpha_j)\\
g(\alpha_i)+g(\alpha_{j+1})&=d(\alpha_i,\alpha_{j+1})
\end{align*}
are equivalent. Hence $g_\varphi=g_\psi$. Analogously for $(i,j+1)\notin
\mathfrak X_N^\circ$ with $j+1-i\geqslant k+3$ it follows that
$g_\psi=g_\varphi$. Hence for $\lambda\in\mathbb Y_N$ we get
$g_\lambda=g_{\lambda^\circ}\in\InjHull(C_N)$.

As a summary, we have the following theorem.
\begin{theorem}
Let $N=2k+1$ be an odd positive integer. The injective hull
\,$\InjHull(C_N)$ of an $N$-cycle  $C_N=\{\alpha_0,\dots,\alpha_{N-1}\}$ with
metric $d(\alpha_i,\alpha_j)=2\min\bigl\{|j-i|,N-|j-i|\bigr\}$ has the
following realization:
The vertices of $\InjHull(C_N)$ are parametrized by the partitions in
$\mathbb Y_N^\circ$:
\begin{align*}
\lambda\in\mathbb Y_N^\circ&\rightsquigarrow\mathcal L_\lambda\subseteq
\mathfrak X_N^\circ&&\mbox{its outer rim
$\mathcal L_\lambda=\bigl\{(i_l,j_l)\bigm|l=0,\dots,N-1\bigr\}$}\\
&\rightsquigarrow g_\lambda\in\mathbb R^{C_N}
&&\mbox{the solution of $g(\alpha_{i_l})+g(\alpha_{j_l})
=d(\alpha_{i_l},\alpha_{j_l})$
for $l=0,\dots,N-1$,}\\
&&&\mbox{explicitly: $g_\lambda(\alpha_j)=\bigl|\tau^j(\lambda)\bigr|-
\frac12k(k-1)$.}
\end{align*}
The $v$-faces of $\InjHull(C_N)$ are the convex hulls of the $2^v$ vertices
$g_\nu$ where $\nu\in[\lambda',\lambda]$ with $\lambda,\lambda'\in
\mathbb Y_N^\circ$ and $\lambda'$ is got by removing
$v$ inner corners from $\lambda$.
\end{theorem}

\begin{remark}
From the explicit description of $\InjHull(C_N)$ in the theorem above
we get
$$\lambda\in\mathbb Y_N^\circ\Longrightarrow\bigl|\tau^j(\lambda)\bigr|
\geqslant\tfrac12 k(k-1)\quad\mbox{for $j=0,\dots,N-1$; and recall $N=2k+1$}.$$
But these inequalities do not characterize the partitions in
$\mathbb Y_N^\circ$ among those in $\mathbb Y_N$, e.\,g.
$$\begin{array}{c||ccccccc}
\multicolumn{8}{l}{
\mbox{for $N=7$ with $\lambda=(3)\notin\mathbb Y_7^\circ$}}\\[2mm]
j\phantom{\Bigl|}&0&1&2&3&4&5&6\\\hline
\mathstrut|\tau^j(3)|-3\phantom{\Bigl|}&0&4&6&6&6&4&0
\end{array}\qquad
\begin{array}{c||ccccccccc}
\multicolumn{10}{l}{
\mbox{for $N=9$ with $\lambda=(4,3)\notin\mathbb Y_9^\circ$}}\\[2mm]
j\phantom{\Bigl|}&0&1&2&3&4&5&6&7&8\\\hline
\mathstrut|\tau^j(4,3)|-6\phantom{\Bigl|}&1&5&7&7&7&7&5&1&1
\end{array}$$
Here is another thought: For $\lambda\in\mathbb Y_N$ let
$$\Vert\lambda\Vert_1:=\sum_{j=0}^{N-1}\bigl|\tau^j(\lambda)\bigr|.$$
The fixed point $(k,\dots,1)$ attains the minimum
$\min\bigl\{\Vert\lambda\Vert_1\bigm|\lambda\in\mathbb Y_N\bigr\}
=\frac18\bigl(N^3-N\bigr)$, and the least element
$\alpha_0=(k-1,\dots,1)\in\mathbb Y_N^\circ$
has $\Vert\alpha_0\Vert_1=\frac18\bigl(N^3+3N-4\bigr)$.
But $\Vert\lambda\Vert_1\leqslant\Vert\alpha_0\Vert_1$ does not
characterize the partitions in $\mathbb Y_N^\circ$, e.\,g.
\begin{align*}
&\mbox{for $N=11$: $\Vert(5,4,3)\Vert_1=170=\Vert\alpha_0\Vert_1$,
but $(5,4,3)\notin\mathbb Y_{11}^\circ$,}\\
&\mbox{for $N=13$: $\Vert(6,5,4,3)\Vert_1=278<279=\Vert\alpha_0\Vert_1$,
but $(6,5,4,3)\notin\mathbb Y_{13}^\circ$.}
\end{align*}
\end{remark}

\section{Census of faces in $\InjHull(C_N)$}
We have to count the homotopically nontrivial loops of length $N$ in
the M\"obius strip $\mathfrak X_N^\circ$. Moreover, for each such loop
$\mathcal L_\lambda$ we need to keep records on the number of 
$\circ$-inner corners of $\lambda$, defined as follows.
\begin{definition}\label{circinner}
An inner corner of $\lambda\in\mathbb Y_N^\circ$ is called \emph{$\circ$-inner}
if its removal from the Young diagram of $\lambda$ results in the Young
diagram of a partition in $\mathbb Y_N^\circ$.
\end{definition}

Rather than trying to extract
results from the literature on lattice paths, we give a self-contained
approach by doing an induction on $k$ in $N=2k+1$, beginning with the base
case $k=1$. We construct a $3$\/$\times$\/$3$ 
matrix $Z_k$ whose entries count certain admissible paths in a strip 
according to their number of $\circ$-inner corners (in an evident sense).
A summand $t^s$ stands for a path with
$s$ \mbox{$\circ$-inner corners}. The first row stands for the
paths passing via $(k-1,N-1)$ and $(k-1,N)$, the second row for those
passing via $(k,N-1)$ and $(k,N)$, and the third row for those
passing via $(k+1,N-1)$ and $(k+1,N)$. These three cases
exhaust all possibilities and are mutually disjoint.
Similarly, the first column stands for the paths passing via $(0,k)$ and
$(1,k)$, the second column for those passing via $(0,k+1)$ and $(1,k+1)$,
and the third column for those passing via $(0,k+2)$ and $(1,k+2)$.
The matrix for the base case is $Z_1=Z$.
\setlength{\unitlength}{0.0003in}
\begin{center}
\begin{picture}(9624,6942)(0,-10)
\texture{55888888 88555555 5522a222 a2555555 55888888 88555555 552a2a2a 2a555555 
	55888888 88555555 55a222a2 22555555 55888888 88555555 552a2a2a 2a555555 
	55888888 88555555 5522a222 a2555555 55888888 88555555 552a2a2a 2a555555 
	55888888 88555555 55a222a2 22555555 55888888 88555555 552a2a2a 2a555555 }
\shade\path(312,6315)(912,6915)(1212,6615)
	(912,6315)(1512,5715)(1212,5415)(312,6315)
\path(312,6315)(912,6915)(1212,6615)
	(912,6315)(1512,5715)(1212,5415)(312,6315)
\shade\path(4212,6015)(4812,6615)(5112,6315)
	(4812,6015)(5112,5715)(4812,5415)(4212,6015)
\path(4212,6015)(4812,6615)(5112,6315)
	(4812,6015)(5112,5715)(4812,5415)(4212,6015)
\shade\path(7812,6015)(8112,6315)(8412,6015)
	(8712,6315)(9012,6015)(8412,5415)(7812,6015)
\path(7812,6015)(8112,6315)(8412,6015)
	(8712,6315)(9012,6015)(8412,5415)(7812,6015)
\shade\path(312,3615)(912,4215)(1812,3315)
	(1512,3015)(912,3615)(612,3315)(312,3615)
\path(312,3615)(912,4215)(1812,3315)
	(1512,3015)(912,3615)(612,3315)(312,3615)
\shade\path(4212,3315)(4812,3915)(5412,3315)
	(5112,3015)(4812,3315)(4512,3015)(4212,3315)
\path(4212,3315)(4812,3915)(5412,3315)
	(5112,3015)(4812,3315)(4512,3015)(4212,3315)
\shade\path(8112,3015)(8412,3315)(8112,3615)
	(8412,3915)(9012,3315)(8412,2715)(8112,3015)
\path(8112,3015)(8412,3315)(8112,3615)
	(8412,3915)(9012,3315)(8412,2715)(8112,3015)
\shade\path(312,915)(912,1515)(1212,1215)
	(1512,1515)(2112,915)(1812,615)
	(1512,915)(1212,615)(912,915)
	(612,615)(312,915)
\path(312,915)(912,1515)(1212,1215)
	(1512,1515)(2112,915)(1812,615)
	(1512,915)(1212,615)(912,915)
	(612,615)(312,915)
\shade\path(4212,615)(5112,1515)(5712,915)
	(5412,615)(5112,915)(4512,315)(4212,615)
\path(4212,615)(5112,1515)(5712,915)
	(5412,615)(5112,915)(4512,315)(4212,615)
\shade\path(8112,315)(8712,915)(8412,1215)
	(8712,1515)(9312,915)(8412,15)(8112,315)
\path(8112,315)(8712,915)(8412,1215)
	(8712,1515)(9312,915)(8412,15)(8112,315)
\path(12,6615)(1212,5415)(2412,6615)
	(2112,6915)(912,5715)
\path(2112,6315)(1512,6915)(612,6015)
\path(1812,6015)(912,6915)(312,6315)
\path(1512,5715)(312,6915)(12,6615)
\path(3612,6615)(4812,5415)(6012,6615)
	(5712,6915)(4512,5715)
\path(5712,6315)(5112,6915)(4212,6015)
\path(5412,6015)(4512,6915)(3912,6315)
\path(5112,5715)(3912,6915)(3612,6615)
\path(7212,6615)(8412,5415)(9612,6615)
	(9312,6915)(8112,5715)
\path(9312,6315)(8712,6915)(7812,6015)
\path(9012,6015)(8112,6915)(7512,6315)
\path(8712,5715)(7512,6915)(7212,6615)
\path(12,3915)(1212,2715)(2412,3915)
	(2112,4215)(912,3015)
\path(2112,3615)(1512,4215)(612,3315)
\path(1812,3315)(912,4215)(312,3615)
\path(1512,3015)(312,4215)(12,3915)
\path(3612,3915)(4812,2715)(6012,3915)
	(5712,4215)(4512,3015)
\path(5712,3615)(5112,4215)(4212,3315)
\path(5412,3315)(4512,4215)(3912,3615)
\path(5112,3015)(3912,4215)(3612,3915)
\path(7212,3915)(8412,2715)(9612,3915)
	(9312,4215)(8112,3015)
\path(9312,3615)(8712,4215)(7812,3315)
\path(9012,3315)(8112,4215)(7512,3615)
\path(8712,3015)(7512,4215)(7212,3915)
\path(12,1215)(1212,15)(2412,1215)
	(2112,1515)(912,315)
\path(2112,915)(1512,1515)(612,615)
\path(1812,615)(912,1515)(312,915)
\path(1512,315)(312,1515)(12,1215)
\path(3612,1215)(4812,15)(6012,1215)
	(5712,1515)(4512,315)
\path(5712,915)(5112,1515)(4212,615)
\path(5412,615)(4512,1515)(3912,915)
\path(5112,315)(3912,1515)(3612,1215)
\path(7212,1215)(8412,15)(9612,1215)
	(9312,1515)(8112,315)
\path(9312,915)(8712,1515)(7812,615)
\path(9012,615)(8112,1515)(7512,915)
\path(8712,315)(7512,1515)(7212,1215)
\put(5712,5415){\makebox(0,0)[lb]{$0$}}
\put(5712,2715){\makebox(0,0)[lb]{$t$}}
\put(5712,15){\makebox(0,0)[lb]{$t$}}
\put(9312,5415){\makebox(0,0)[lb]{$1$}}
\put(9312,2715){\makebox(0,0)[lb]{$0$}}
\put(9312,15){\makebox(0,0)[lb]{$0$}}
\put(2112,5415){\makebox(0,0)[lb]{$0$}}
\put(2112,2715){\makebox(0,0)[lb]{$t$}}
\put(2112,15){\makebox(0,0)[lb]{$t^2$}}
\end{picture}\qquad\qquad\raisebox{24mm}[0pt][0pt]{$Z=\begin{pmatrix}
0&0&1\\t&t&0\\t^2&t&0
\end{pmatrix}$}
\end{center}
The matrix $S$ incorporates the induction from $k$ to $k+1$. 
The first row stands for the paths passing via $(k,N)$ and $(k,N+1)$,
the second row for those passing via $(k+1,N)$ and $(k+1,N+1)$, and the
third row for those passing via $(k+2,N)$ and $(k+2,N+1)$.
Similarly, the first column stands for the paths passing via
$(k-1,N-1)$ and $(k-1,N)$, the second column for those
passing via $(k,N-1)$ and $(k,N)$, and the third column for those
passing via $(k+1,N-1)$ and $(k+1,N)$. (Note that to get the site positions
in $\mathfrak X_{N+2}^\circ\subseteq\mathfrak X_{N+2}$, we have to shift by
$(0,1)$.)
\begin{center}
\begin{picture}(9624,6942)(0,-10)
\texture{55888888 88555555 5522a222 a2555555 55888888 88555555 552a2a2a 2a555555 
	55888888 88555555 55a222a2 22555555 55888888 88555555 552a2a2a 2a555555 
	55888888 88555555 5522a222 a2555555 55888888 88555555 552a2a2a 2a555555 
	55888888 88555555 55a222a2 22555555 55888888 88555555 552a2a2a 2a555555 }
\shade\path(12,3315)(312,3615)(612,3315)
	(1212,3915)(1812,3315)(1512,3015)
	(1212,3315)(612,2715)(12,3315)
\path(12,3315)(312,3615)(612,3315)
	(1212,3915)(1812,3315)(1512,3015)
	(1212,3315)(612,2715)(12,3315)
\shade\path(12,6015)(312,6315)(612,6015)
	(912,6315)(1512,5715)(1212,5415)
	(912,5715)(612,5415)(12,6015)
\path(12,6015)(312,6315)(612,6015)
	(912,6315)(1512,5715)(1212,5415)
	(912,5715)(612,5415)(12,6015)
\shade\path(3912,6315)(4212,6615)(5112,5715)
	(4812,5415)(3912,6315)
\path(3912,6315)(4212,6615)(5112,5715)
	(4812,5415)(3912,6315)
\shade\path(7812,6615)(8112,6915)(8712,6315)
	(8412,6015)(8712,5715)(8412,5415)
	(7812,6015)(8112,6315)(7812,6615)
\path(7812,6615)(8112,6915)(8712,6315)
	(8412,6015)(8712,5715)(8412,5415)
	(7812,6015)(8112,6315)(7812,6615)
\shade\path(3912,3615)(4212,3915)(4512,3615)
	(4812,3915)(5412,3315)(5112,3015)
	(4812,3315)(4512,3015)(3912,3615)
\path(3912,3615)(4212,3915)(4512,3615)
	(4812,3915)(5412,3315)(5112,3015)
	(4812,3315)(4512,3015)(3912,3615)
\shade\path(7812,3915)(8112,4215)(9012,3315)
	(8712,3015)(7812,3915)
\path(7812,3915)(8112,4215)(9012,3315)
	(8712,3015)(7812,3915)
\shade\path(12,615)(312,915)(612,615)
	(1512,1515)(2112,915)(1812,615)
	(1512,915)(612,15)(12,615)
\path(12,615)(312,915)(612,615)
	(1512,1515)(2112,915)(1812,615)
	(1512,915)(612,15)(12,615)
\shade\path(3912,915)(4212,1215)(4512,915)
	(5112,1515)(5712,915)(5412,615)
	(5112,915)(4512,315)(3912,915)
\path(3912,915)(4212,1215)(4512,915)
	(5112,1515)(5712,915)(5412,615)
	(5112,915)(4512,315)(3912,915)
\shade\path(7812,1215)(8112,1515)(8412,1215)
	(8712,1515)(9312,915)(9012,615)
	(8712,915)(8412,615)(7812,1215)
\path(7812,1215)(8112,1515)(8412,1215)
	(8712,1515)(9312,915)(9012,615)
	(8712,915)(8412,615)(7812,1215)
\path(1512,5715)(612,6615)
\path(2112,6315)(1512,6915)(312,5715)
\path(1812,6015)(912,6915)(12,6015)
	(612,5415)(2112,6915)(2412,6615)
	(1212,5415)(312,6315)
\path(5112,5715)(4212,6615)
\path(5712,6315)(5112,6915)(3912,5715)
\path(5412,6015)(4512,6915)(3612,6015)
	(4212,5415)(5712,6915)(6012,6615)
	(4812,5415)(3912,6315)
\path(8712,5715)(7812,6615)
\path(9312,6315)(8712,6915)(7512,5715)
\path(9012,6015)(8112,6915)(7212,6015)
	(7812,5415)(9312,6915)(9612,6615)
	(8412,5415)(7512,6315)
\path(1512,3015)(612,3915)
\path(2112,3615)(1512,4215)(312,3015)
\path(1812,3315)(912,4215)(12,3315)
	(612,2715)(2112,4215)(2412,3915)
	(1212,2715)(312,3615)
\path(1512,315)(612,1215)
\path(2112,915)(1512,1515)(312,315)
\path(1812,615)(912,1515)(12,615)
	(612,15)(2112,1515)(2412,1215)
	(1212,15)(312,915)
\path(5112,3015)(4212,3915)
\path(5712,3615)(5112,4215)(3912,3015)
\path(5412,3315)(4512,4215)(3612,3315)
	(4212,2715)(5712,4215)(6012,3915)
	(4812,2715)(3912,3615)
\path(5112,315)(4212,1215)
\path(5712,915)(5112,1515)(3912,315)
\path(5412,615)(4512,1515)(3612,615)
	(4212,15)(5712,1515)(6012,1215)
	(4812,15)(3912,915)
\path(8712,3015)(7812,3915)
\path(9312,3615)(8712,4215)(7512,3015)
\path(9012,3315)(8112,4215)(7212,3315)
	(7812,2715)(9312,4215)(9612,3915)
	(8412,2715)(7512,3615)
\path(8712,315)(7812,1215)
\path(9312,915)(8712,1515)(7512,315)
\path(9012,615)(8112,1515)(7212,615)
	(7812,15)(9312,1515)(9612,1215)
	(8412,15)(7512,915)
\put(2112,5415){\makebox(0,0)[lb]{$1$}}
\put(5712,5415){\makebox(0,0)[lb]{$1$}}
\put(9312,5415){\makebox(0,0)[lb]{$0$}}
\put(2112,2715){\makebox(0,0)[lb]{$t$}}
\put(5712,2715){\makebox(0,0)[lb]{$t$}}
\put(9312,2715){\makebox(0,0)[lb]{$1$}}
\put(2112,15){\makebox(0,0)[lb]{$t$}}
\put(5712,15){\makebox(0,0)[lb]{$t$}}
\put(9312,15){\makebox(0,0)[lb]{$t$}}
\end{picture}\qquad\qquad\raisebox{24mm}[0pt][0pt]{$S=\begin{pmatrix}
1&1&0\\t&t&1\\t&t&t
\end{pmatrix}$}
\end{center}
So $Z_k=S^{k-1}Z$.
Only the paths that pass via $(0,k)$ and $(1,k)$ and $(k+1,N-1)$ (and
$(k+1,N)$) cannot be completed to an admissible loop in the M\"obius
strip $\mathfrak X_N^\circ$. The enumerator for all admissible loops is
therefore the sum
of all matrix entries of $Z_k$ except the entry in the third row and
first column. In other words,
with
$$A=\begin{pmatrix}
1&1&0\\1&1&1\\1&1&1
\end{pmatrix}$$
we obtain
\begin{align}\label{enumerator}
\sum_{\lambda\in\mathbb Y_N^\circ}t^{\#(\textup{$\circ$-inner corners of $\lambda$})}
=\tr(AZ_k)=\tr(AS^{k-1}Z)=\tr(S^{k-1}ZA).
\end{align}

\begin{lemma}\label{Spowers}
For $n\geqslant1$ the matrix powers of $S$ are
$$S^n=a_n S+b_n (S^2-(1+t)S)$$
where
\begin{align*}
a_n&=\sum_{j=0}^{n-1}\binom{2(n-1)-j}{j}t^j\quad\mbox{and}\quad
b_n=\sum_{j=0}^{n-2}\binom{2(n-1)-j-1}{j}t^j.
\end{align*}
\end{lemma}
\begin{prf}
We do an induction on $n$. The formula is true for $n=1$. For $n\geqslant1$ we
compute using $S^3=(1+2t)S^2-t^2S$
\begin{align*}
S^{n+1}&=S\bigl(a_n S+b_n (S^2-(1+t)S)\bigr)
=a_nS^2+b_n((1+2t)S^2-t^2S-(1+t)S^2)\\
&=((1+t)a_n+tb_n)S+(a_n+tb_n)(S^2-(1+t)S)
\end{align*}
and hence we must show that
\begin{align}\label{abrecursion}
a_{n+1}&=(1+t)a_n+tb_n\quad\mbox{and}\quad
b_{n+1}=a_n+tb_n.
\end{align}
The constant coefficients of the four expressions in (\ref{abrecursion}) are
$1$, and for the coefficient of $t^{j+1}$ we have
\begin{align}\notag
[t^{j+1}](a_n+tb_n)&=\binom{2(n-1)-(j+1)}{j+1}+\binom{2(n-1)-j-1}{j}\\
\label{bn+1}
&=\binom{2n-(j+1)-1}{j+1}=[t^{j+1}]b_{n+1}
\intertext{and using (\ref{bn+1}) we continue}\notag
[t^{j+1}]((1+t)a_n+tb_n)&=[t^{j+1}](b_{n+1}+ta_n)\\\notag
&=\binom{2n-(j+1)-1}{j+1}+\binom{2(n-1)-j}{j}\\\notag
&=\binom{2n-(j+1)}{j+1}=[t^{j+1}]a_{n+1}
\end{align}
and hence (\ref{abrecursion}) is verified.
\end{prf}

We continue with the computation in (\ref{enumerator})
using $ZA=S^2-tS$ and Lemma~\ref{Spowers}.
\begin{align}\label{ASkZ}
\raisebox{0pt}[0pt][0pt]{$\displaystyle
\sum_{\lambda\in\mathbb Y_N^\circ}t^{\#(\textup{$\circ$-inner corners of $\lambda$})}$}
&=\tr\bigl(S^{k-1}ZA\bigr)=
\tr\bigl(S^{k+1}-tS^k\bigr)\\\notag
&=\tr\bigl(a_{k+1}S+b_{k+1}(S^2-(1+t)S)-ta_kS-tb_k(S^2-(1+t)S)\bigr)\\\notag
&=(a_{k+1}-ta_k)\tr(S)+(b_{k+1}-tb_k)\tr(S^2-(1+t)S)\\\notag
&=(a_{k+1}-ta_k)(1+2t)+(b_{k+1}-tb_k)t\\\notag
&=(1+t)(a_{k+1}-ta_k)+t(b_{k+1}-tb_k)+t(a_{k+1}-ta_k)
\intertext{and according to (\ref{abrecursion}) this can be written as}\notag
&=a_{k+2}-t^2a_k=:U.
\end{align}
The coefficients of $t^s$ are $[t^0]U=1$, $[t^1]U=2k+1$ and for
$s\geqslant2$
\begin{align}\notag
[t^s]U&=\binom{2(k+1)-s}{s}-\binom{2(k-1)-(s-2)}{s-2}
\stackrel{(\ast)}{=}\binom{2k-s}{s-1}+\binom{2k+1-s}{s}\\\label{formU}
&=\Bigl(\frac{s}{2k+1-s}+1\Bigr)\binom{2k+1-s}{s}
=\frac{2k+1}{2k+1-s}\binom{2k+1-s}{s}
\end{align}
where $(\ast)$ is quickly seen from Pascal's triangle, symbolically
$$\left\{\begin{array}{c@{}c@{}c@{}c}
\phantom{+}&\phantom{+}&\phantom{+}&\phantom{+}\\[-6mm]
\mathstrut-&&\cdot\\&\mathstrut\cdot&&\cdot\\&&\mathstrut+
\end{array}\right\}
\quad=\quad\left\{\begin{array}{c@{}c@{}c@{}c}
\phantom{+}&\phantom{+}&\phantom{+}&\phantom{+}\\[-6mm]
\mathstrut-&&\cdot\\&\mathstrut+&&+\\\mathstrut&&\cdot
\end{array}\right\}
\quad=\quad\left\{\begin{array}{c@{}c@{}c@{}c}
\phantom{+}&\phantom{+}&\phantom{+}&\phantom{+}\\[-6mm]
\mathstrut\cdot&&+\\&\mathstrut\cdot&&\mathstrut+\\&&\cdot
\end{array}\right\}$$
and since (\ref{formU}) is also the correct expression for $s=0$ and $s=1$,
we have
\begin{align}\label{admcornerenumerator}
\sum_{\lambda\in\mathbb Y_N^\circ}t^{\#(\textup{$\circ$-inner corners of $\lambda$})}
=\sum_{s=0}^{(N-1)/2}
\frac{N}{N-s}\binom{N-s}{s}t^s.
\end{align}

As an alternative, we compute
\begin{align*}
(\mathbf{1}-qS)^{-1}=\frac{1}{1-(1+2t)q+t^2q^2}
\begin{pmatrix}
1-2tq-(t-t^2)q^2&q-tq^2&q^2\\
tq+(t-t^2)q^2&\makebox[28mm]{$1-(1+t)q+tq^2$}&q-q^2\\
tq&tq&1-(1+t)q
\end{pmatrix}
\end{align*}
from which we get the generating series
$$\sum_{k=0}^\infty\tr(S^k)\,q^k=\frac{3-(2+4t)q+t^2q^2}{1-(1+2t)q+t^2q^2}=
3+\frac{(1+2t)q-2t^2q^2}{1-(1+2t)q+t^2q^2}=
3+(1+2t)q+\cdots$$
and hence
\begin{align*}
1+\sum_{k=1}^\infty\tr\bigl(S^{k+1}\!-\!tS^k\bigr)\,q^k&=
1+\Bigl(\frac{(1+2t)-2t^2q}{1-(1\!+\!2t)q+t^2q^2}-(1\!+\!2t)\Bigr)
-t\frac{(1+2t)q-2t^2q^2}{1-(1\!+\!2t)q+t^2q^2}\\
&=\frac{1+tq}{1-(1+2t)q+t^2q^2}
\end{align*}
\begin{align}\notag
\makebox[4cm]{}&=\frac{\dfrac{1+\sqrt{1+4t}}{2}}{1-
\left(\dfrac{1+\sqrt{1+4t}}{2}\right)^{2\mathstrut}q}+
\frac{\dfrac{1-\sqrt{1+4t}}{2}}{1-
\left(\dfrac{1-\sqrt{1+4t}}{2}\right)^{2\mathstrut}q}\\[5mm]\label{genseriesS}
&=\sum_{k=0}^\infty\left(\left(\frac{1+\sqrt{1+4t}}{2}\right)^{2k+1}
+\left(\frac{1-\sqrt{1+4t}}{2}\right)^{2k+1}\right)q^k.
\end{align}
From (\ref{ASkZ}) and the
generating series (\ref{genseriesS}) we have (at least for $N\geqslant3$)
\begin{align}\label{enumerator2}
\sum_{\lambda\in\mathbb Y_N^\circ}t^{\#(\textup{$\circ$-inner corners of $\lambda$})}
=\left(\frac{1+\sqrt{1+4t}}{2}\right)^N
+\left(\frac{1-\sqrt{1+4t}}{2}\right)^N.
\end{align}

\begin{theorem}
Let $N=2k+1$ be an odd positive integer. The injective hull \,$\InjHull(C_N)$
of an $N$-cycle satisfies
\begin{align}\label{vfacesgs1}
\sum_{v\geqslant0}\mbox{\textup{\#($v$-faces in $\InjHull(C_N)$)}}\,t^v&=
\left(\frac{1+\sqrt{5+4t}}{2}\right)^N+
\left(\frac{1-\sqrt{5+4t}}{2}\right)^N\\\label{vfacesgs2}
&=\sum_{s=0}^{(N-1)/2}\frac{N}{N-s}\binom{N-s}{s}(1+t)^s
\intertext{and it follows that}\label{vfacesformula}
\mbox{\textup{\#($v$-faces in $\InjHull(C_N)$)}}&=
\frac{1}{2^{N-2v-1}}\sum_{s=v}^{(N-1)/2}\binom{N}{2s}\binom{s}{v}\,5^{s-v}\\\label{vfacesformula2}
&=\sum_{s=v}^{(N-1)/2}\frac{N}{N-s}\binom{N-s}{s}\binom{s}{v}.
\end{align}
\end{theorem}
\begin{prf}
$\InjHull(C_1)$ consists of a single point, in accordance with the
statement of the theorem. Let us suppose that $k\geqslant1$.
Recall that each $\lambda\in\mathbb Y_N^\circ$ with $s$ $\circ$-inner
corners gives rise to $\binom sv$ $v$-faces with ``top vertex''
$g_\lambda\in\InjHull(C_N)$. Hence 
$$\sum_{v\geqslant0}\mbox{\textup{\#($v$-faces in $\InjHull(C_N)$)}}\,t^v
=\sum_{\lambda\in\mathbb Y_N^\circ}(1+t)^{\#(\textup{$\circ$-inner corners
of $\lambda$})}$$
and (\ref{vfacesgs1}) follows from (\ref{enumerator2}), whereas
(\ref{vfacesgs2}) follows from (\ref{admcornerenumerator}).
We expand the right hand side of (\ref{vfacesgs1})
\begin{align*}
\left(\frac{1+\sqrt{5+4t}}{2}\right)^N+
\left(\frac{1-\sqrt{5+4t}}{2}\right)^N&=\frac{1}{2^{N-1}}\sum_{s=0}^{(N-1)/2}
\binom{N}{2s}\sum_{v=0}^s\binom sv 5^{s-v}(4t)^v
\end{align*}
from which we have (\ref{vfacesformula}). Finally,
(\ref{vfacesformula2}) follows directly from (\ref{vfacesgs2}).
\end{prf}

\begin{remark}\label{LucasL}
The number of vertices in $\InjHull(C_N)$, still for $N$ odd, is the $N$th
Lucas number
$L_N=\varphi^N+(1-\varphi)^N$, where $\varphi=\frac12(1+\sqrt5)$ is the
golden ratio, which follows by putting $t=0$ in (\ref{vfacesgs1}). Since the
cyclic group of order $N$ permutes the vertices of $\InjHull(C_N)$ with
exactly one fixed point, one has the congruence $L_p\equiv1\pmod p$
for each odd prime number $p$, a result that is also evident by using
either of the expressions (\ref{vfacesformula}) and (\ref{vfacesformula2}) 
specialized to $v=0$ for the odd-indexed Lucas numbers.

The total number of (nonempty) faces in $\InjHull(C_N)$ is $2^N-1$, which
follows by putting $t=1$ in (\ref{vfacesgs1}). Another special value to
insert in (\ref{vfacesgs1}) is $t=-1$. The Euler characteristic of
$\InjHull(C_N)$ is $1$, as it must be because injective hulls are
contractible.
\end{remark}

\begin{remark}
Let us briefly give the analogous computations for $N=2k$ even.
The discrete M\"obius strip is $\mathfrak X_N^\circ:=\bigl\{(i,j)\in
\mathfrak X_N\bigm|d(\alpha_i,\alpha_j)\geqslant k-1\bigr\}$ and
has $\frac32 N$ sites. Let $\mathbb Y_N^\circ$ be defined as in
Definition~\ref{YNcirc}. The matrices that are used for the census are
\setlength{\unitlength}{0.0003in}
\begin{center}
\begin{picture}(5534,3942)(0,-10)
\texture{55888888 88555555 5522a222 a2555555 55888888 88555555 552a2a2a 2a555555 
	55888888 88555555 55a222a2 22555555 55888888 88555555 552a2a2a 2a555555 
	55888888 88555555 5522a222 a2555555 55888888 88555555 552a2a2a 2a555555 
	55888888 88555555 55a222a2 22555555 55888888 88555555 552a2a2a 2a555555 }
\shade\path(312,3315)(912,3915)(1212,3615)
	(912,3315)(1212,3015)(912,2715)(312,3315)
\path(312,3315)(912,3915)(1212,3615)
	(912,3315)(1212,3015)(912,2715)(312,3315)
\shade\path(3912,3315)(4212,3615)(4512,3315)
	(4812,3615)(5112,3315)(4512,2715)(3912,3315)
\path(3912,3315)(4212,3615)(4512,3315)
	(4812,3615)(5112,3315)(4512,2715)(3912,3315)
\shade\path(312,615)(912,1215)(1512,615)
	(1212,315)(912,615)(612,315)(312,615)
\path(312,615)(912,1215)(1512,615)
	(1212,315)(912,615)(612,315)(312,615)
\shade\path(4212,315)(4512,615)(4212,915)
	(4512,1215)(5112,615)(4512,15)(4212,315)
\path(4212,315)(4512,615)(4212,915)
	(4512,1215)(5112,615)(4512,15)(4212,315)
\path(12,3615)(912,2715)(1812,3615)
	(1512,3915)(612,3015)
\path(1512,3315)(912,3915)(312,3315)
\path(1212,3015)(312,3915)(12,3615)
\path(3612,3615)(4512,2715)(5412,3615)
	(5112,3915)(4212,3015)
\path(5112,3315)(4512,3915)(3912,3315)
\path(4812,3015)(3912,3915)(3612,3615)
\path(12,915)(912,15)(1812,915)
	(1512,1215)(612,315)
\path(1512,615)(912,1215)(312,615)
\path(1212,315)(312,1215)(12,915)
\path(3612,915)(4512,15)(5412,915)
	(5112,1215)(4212,315)
\path(5112,615)(4512,1215)(3912,615)
\path(4812,315)(3912,1215)(3612,915)
\put(6000,1865){\makebox(0,0)[l]{$Z=\begin{pmatrix}0&1\\t&0\end{pmatrix}$}}
\end{picture}
\makebox[4cm]{}
\begin{picture}(8834,3942)(0,-10)
\texture{55888888 88555555 5522a222 a2555555 55888888 88555555 552a2a2a 2a555555 
	55888888 88555555 55a222a2 22555555 55888888 88555555 552a2a2a 2a555555 
	55888888 88555555 5522a222 a2555555 55888888 88555555 552a2a2a 2a555555 
	55888888 88555555 55a222a2 22555555 55888888 88555555 552a2a2a 2a555555 }
\shade\path(12,615)(312,915)(612,615)
	(1212,1215)(1812,615)(1512,315)
	(1212,615)(612,15)(12,615)
\path(12,615)(312,915)(612,615)
	(1212,1215)(1812,615)(1512,315)
	(1212,615)(612,15)(12,615)
\shade\path(12,3315)(312,3615)(612,3315)
	(912,3615)(1512,3015)(1212,2715)
	(912,3015)(612,2715)(12,3315)
\path(12,3315)(312,3615)(612,3315)
	(912,3615)(1512,3015)(1212,2715)
	(912,3015)(612,2715)(12,3315)
\shade\path(3912,3615)(4212,3915)(5112,3015)
	(4812,2715)(3912,3615)
\path(3912,3615)(4212,3915)(5112,3015)
	(4812,2715)(3912,3615)
\shade\path(3912,915)(4212,1215)(4512,915)
	(4812,1215)(5412,615)(5112,315)
	(4812,615)(4512,315)(3912,915)
\path(3912,915)(4212,1215)(4512,915)
	(4812,1215)(5412,615)(5112,315)
	(4812,615)(4512,315)(3912,915)
\path(12,3315)(612,2715)(1812,3915)
	(2112,3615)(1212,2715)(312,3615)(12,3315)
\path(1812,3315)(1212,3915)(312,3015)
\path(1512,3015)(612,3915)(12,3315)
\path(3612,3315)(4212,2715)(5412,3915)
	(5712,3615)(4812,2715)(3912,3615)(3612,3315)
\path(5412,3315)(4812,3915)(3912,3015)
\path(12,615)(612,15)(1812,1215)
	(2112,915)(1212,15)(312,915)(12,615)
\path(1512,315)(612,1215)(12,615)
\path(312,315)(1212,1215)
\path(3612,615)(4212,15)(5412,1215)
	(5712,915)(4812,15)(3912,915)(3612,615)
\path(5112,315)(4512,915)(3912,315)
\put(6000,1865){\makebox(0,0)[l]{$S=\begin{pmatrix}1&1\\t&t\end{pmatrix}$}}
\end{picture}
\end{center}
and $A$ is the $2$\/$\times$\/$2$ matrix with all entries $1$.
Then
$$\sum_{\lambda\in\mathbb Y_N^\circ}t^{\#(\textup{$\circ$-inner corners of
$\lambda$})}
=\tr(AS^{k-1}Z)=\tr(S^{k-1}ZA)=\tr(S^k)=(1+t)^k$$
and hence
$$\sum_{v\geqslant0}\mbox{\textup{\#($v$-faces in $\InjHull(C_{2k})$)}}\,t^v
=(2+t)^k$$
as was clear before because $\InjHull(C_{2k})$ is a $k$-dimensional cube.
\end{remark}

\section{Pictures and additional material}
For $N=2k+1\geqslant5$ we have $\alpha:=\alpha_0=(k-1,\dots,1)
\in\mathbb Y_N^\circ$, which is the least element. 
Let $\beta=(k,\dots,2)\in\mathbb Y_N^\circ$. All the
$2^{k-1}$ partitions in the interval $[\alpha,\beta]$ have
length $k-1$. For such a partition $\lambda$ let $(\lambda\,1)$ denote the
partition with one additional part $1$. In particular,
$(\beta\,1)=(k,\dots,1)\in\mathbb Y_N^\circ$ is
fixed by the cyclic action $\tau$. We have the evident decomposition
$$\bigl[\alpha,(\beta\,1)\bigr]=[\alpha,\beta]\cup\bigl\{(\lambda\,1)\bigm|
\lambda\in[\alpha,\beta]\bigr\}=[\alpha,\beta]\cup
\bigl[(\alpha\,1),(\beta\,1)\bigr].$$ 
The convex hull of the vertices $g_\nu$ for
$\nu\in\bigl[\alpha,(\beta\,1)\bigr]$ is one of the $N$ faces of
$\InjHull(C_N)$ of maximal dimension. The $1$-skeleton of this $k$-cube
is drawn in red in the following pictures for $N=5,7,9,11$.
The part for $\lambda\in[\alpha,\beta]$ is drawn with somewhat thicker
edges.

For $\lambda=\bigl(k-1+\delta_1,\dots,1+\delta_{k-1}\bigr)
\in[\alpha,\beta]$ (i.\,e.\ with $\delta_j\in\{0,1\}$),
which we also write as $\lambda=\alpha\stackrel{.}{+}
(\delta_1,\dots,\delta_{k-1})$, we get
$$\tau(\lambda\,1)=\begin{cases}
\bigl(k,k-2+\delta_1,\dots,1+\delta_{k-2}\bigr)
&\mbox{if $\delta_{k-1}=0$,}\\
\bigl(k,k-2+\delta_1,\dots,1+\delta_{k-2},1\bigr)
&\mbox{if $\delta_{k-1}=1$.}
\end{cases}$$
and rewrite it as
$$\tau\bigl(\alpha\stackrel{.}{+}(\delta_1,\dots,\delta_{k-1})\,1
\bigr)
=\begin{cases}
\alpha\stackrel{.}{+}(1,\delta_1,\dots,\delta_{k-2})
&\mbox{if $\delta_{k-1}=0$,}\\
\bigl(\alpha\stackrel{.}{+}(1,\delta_1,\dots,
\delta_{k-2})\,1\bigr)
&\mbox{if $\delta_{k-1}=1$.}
\end{cases}$$
Hence for each partition $(\lambda\,1)$ with $\lambda\in[\alpha,\beta)$
there is an exponent $l\geqslant0$ such that
$\tau^j(\lambda\,1)=(\lambda_j\,1)$ with $\lambda_j\in[\alpha,\beta)$ for
$j=0,\dots,l$ and $\tau^{l+1}(\lambda\,1)=\lambda_{l+1}\in[\alpha,\beta]$.
In fact, for $\lambda=\alpha\stackrel{.}{+}(\delta_1,\dots,
\delta_{k-1})$ the corresponding exponent is
$l=\min\bigl\{j\bigm|\delta_{k-1-j}=0\bigr\}$.
So every vertex of $\InjHull(C_N)$ that belongs to at
least one of the $N$ $k$-cubes is either the
fixed point $g_{(\beta\,1)}$ or belongs to the orbit of one of the points
$g_\lambda$ for $\lambda\in[\alpha,\beta]$. 

For $\lambda=\alpha\stackrel{.}{+}(\delta_1,\dots,\delta_{k-1})
\in[\alpha,\beta]$ we have
$(1,k),(2,k+1),\dots,(k+1,N-1)\notin\mathcal L_\lambda$ and
$(0,k+2)\in\mathcal L_\lambda$. The following picture
illustrates the situation for $N=9$. We see three copies of a fundamental
domain in the universal covering.
\setlength{\unitlength}{0.0006in}
\begin{center}
\begin{picture}(9044,2749)(0,-10)
\path(322,2712)(22,2412)(1222,1212)
	(2722,2712)(4222,1212)(5722,2712)
	(7222,1212)(8722,2712)(9022,2412)
	(7822,1212)(6322,2712)(4822,1212)
	(3322,2712)(1822,1212)(322,2712)
\path(922,1512)(2122,2712)(3622,1212)
	(5122,2712)(6622,1212)(8122,2712)(8722,2112)
\path(622,1812)(1522,2712)(3022,1212)
	(4522,2712)(6022,1212)(7522,2712)(8422,1812)
\path(322,2112)(922,2712)(2422,1212)
	(3922,2712)(5422,1212)(6922,2712)(8122,1512)
\texture{55888888 88555555 5522a222 a2555555 55888888 88555555 552a2a2a 2a555555 
	55888888 88555555 55a222a2 22555555 55888888 88555555 552a2a2a 2a555555 
	55888888 88555555 5522a222 a2555555 55888888 88555555 552a2a2a 2a555555 
	55888888 88555555 55a222a2 22555555 55888888 88555555 552a2a2a 2a555555 }
\shade\path(322,2112)(622,2412)(1222,1812)
	(1522,2112)(1822,1812)(2122,2112)
	(2422,1812)(2722,2112)(3022,1812)
	(2722,1512)(2422,1812)(2122,1512)
	(1822,1812)(1222,1212)(322,2112)
\path(322,2112)(622,2412)(1222,1812)
	(1522,2112)(1822,1812)(2122,2112)
	(2422,1812)(2722,2112)(3022,1812)
	(2722,1512)(2422,1812)(2122,1512)
	(1822,1812)(1222,1212)(322,2112)
\shade\path(5722,2112)(6022,2412)(6622,1812)
	(6922,2112)(7222,1812)(7522,2112)
	(7822,1812)(8122,2112)(8422,1812)
	(8122,1512)(7822,1812)(7522,1512)
	(7222,1812)(6622,1212)(5722,2112)
\path(5722,2112)(6022,2412)(6622,1812)
	(6922,2112)(7222,1812)(7522,2112)
	(7822,1812)(8122,2112)(8422,1812)
	(8122,1512)(7822,1812)(7522,1512)
	(7222,1812)(6622,1212)(5722,2112)
\shade\path(3022,1812)(3322,1512)(3922,2112)
	(4222,1812)(4522,2112)(4822,1812)
	(5122,2112)(5422,1812)(5722,2112)
	(5422,2412)(5122,2112)(4822,2412)
	(4522,2112)(3922,2712)(3022,1812)
\path(3022,1812)(3322,1512)(3922,2112)
	(4222,1812)(4522,2112)(4822,1812)
	(5122,2112)(5422,1812)(5722,2112)
	(5422,2412)(5122,2112)(4822,2412)
	(4522,2112)(3922,2712)(3022,1812)
\path(1522,612)(1822,312)(1522,12)
	(1222,312)(1522,612)
\shade\path(3022,612)(3322,312)(3022,12)
	(2722,312)(3022,612)
\path(3022,612)(3322,312)(3022,12)
	(2722,312)(3022,612)
\path(7522,612)(7822,312)(7522,12)
	(7222,312)(7522,612)
\path(6022,612)(6322,312)(6022,12)(5722,312)(6022,612)
\put(1522,312){\makebox(0,0){\footnotesize$1$}}
\put(3022,2112){\makebox(0,0){\footnotesize$\delta_1$}}
\put(3022,1512){\makebox(0,0){\footnotesize$\bar\delta_1$}}
\put(322,2412){\makebox(0,0){\footnotesize$\bar\delta_1$}}
\put(5722,1812){\makebox(0,0){\footnotesize$\delta_1$}}
\put(5722,2412){\makebox(0,0){\footnotesize$\bar\delta_1$}}
\put(8422,2112){\makebox(0,0){\footnotesize$\delta_1$}}
\put(2422,2112){\makebox(0,0){\footnotesize$\delta_2$}}
\put(2422,1512){\makebox(0,0){\footnotesize$\bar\delta_2$}}
\put(5122,1812){\makebox(0,0){\footnotesize$\delta_2$}}
\put(5122,2412){\makebox(0,0){\footnotesize$\bar\delta_2$}}
\put(7822,2112){\makebox(0,0){\footnotesize$\delta_2$}}
\put(7822,1512){\makebox(0,0){\footnotesize$\bar\delta_2$}}
\put(1822,2112){\makebox(0,0){\footnotesize$\delta_3$}}
\put(1822,1512){\makebox(0,0){\footnotesize$\bar\delta_3$}}
\put(4522,1812){\makebox(0,0){\footnotesize$\delta_3$}}
\put(4522,2412){\makebox(0,0){\footnotesize$\bar\delta_3$}}
\put(7222,2112){\makebox(0,0){\footnotesize$\delta_3$}}
\put(7222,1512){\makebox(0,0){\footnotesize$\bar\delta_3$}}
\put(2272,312){\makebox(0,0){\footnotesize$=$}}
\put(6022,312){\makebox(0,0){\footnotesize$0$}}
\put(4522,312){\makebox(0,0){\footnotesize$\bar\delta_j=1-\delta_j$}}
\put(6772,312){\makebox(0,0){\footnotesize$=$}}
\end{picture}
\end{center}
The cyclic action $\tau$ is given by
translating the outer rim in the M\"obius strip. Suppose now that we have an
equality $\lambda=\tau^j(\mu)$ with $j\in\{0,\dots,N-1\}$ and
$\lambda,\mu\in[\alpha,\beta]$. From the considerations above it is
clear that such an equality can only hold if $j=0$ and $\lambda=\mu$.
Thus the subcomplex induced from the faces of maximal dimension has
$1+(2k+1)\cdot2^{k-1}$ vertices. For $k=2$ and $k=3$ there are no
further
vertices in $\InjHull(C_{2k+1})$.
\begin{center}
\raisebox{4mm}{\includegraphics[width=5cm]{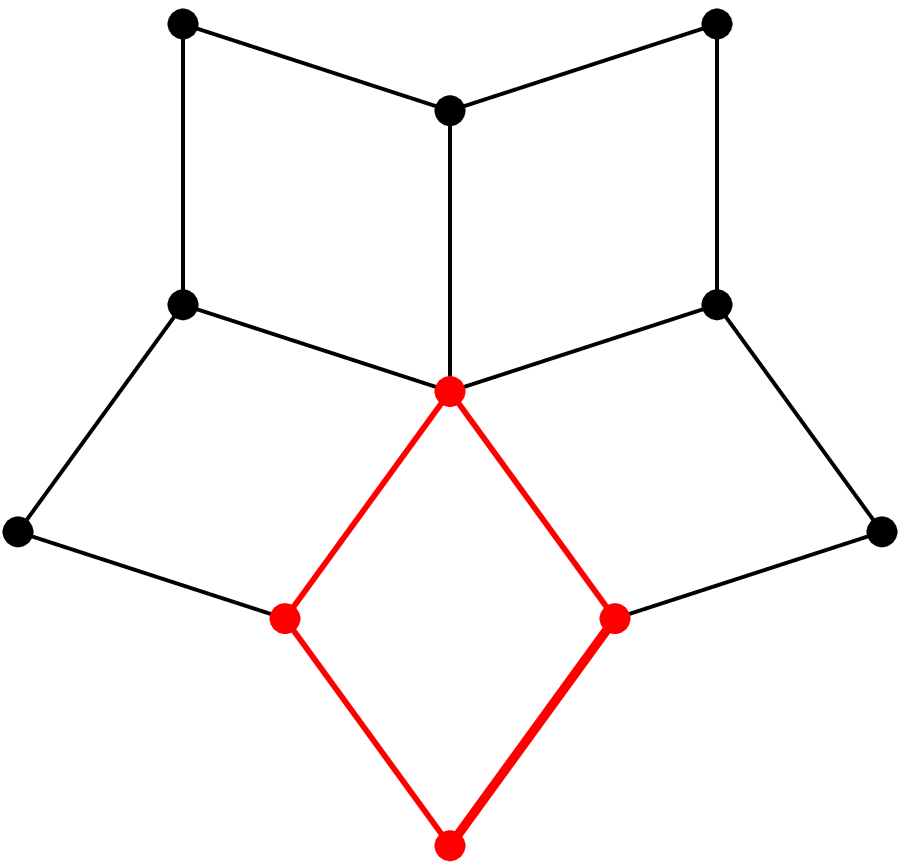}}
\qquad\qquad\includegraphics[width=6cm]{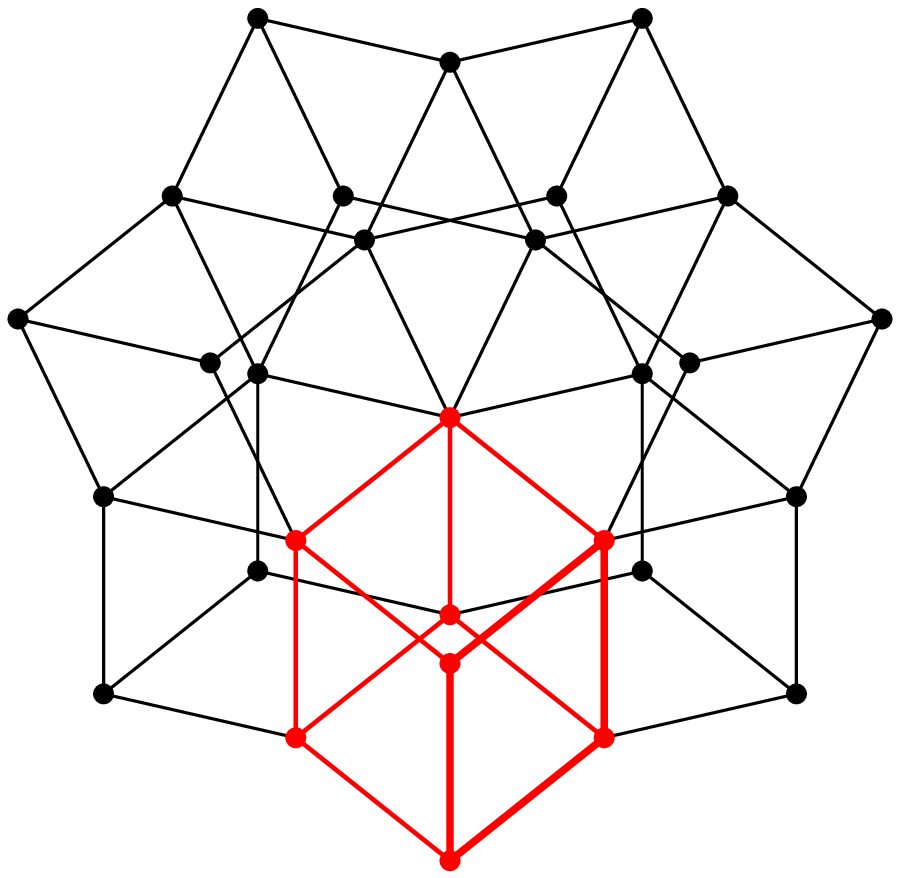}
\end{center}
But the ratio $\bigl(1+(2k+1)\cdot2^{k-1}\bigr)\bigm/L_{2k+1}$ (recall from
Remark~\ref{LucasL} that the total number of vertices in $\InjHull(C_{2k+1})$
is $L_{2k+1}$) tends to zero as $k$ tends to infinity.
The three additional vertices in $\InjHull(C_9)$, that do not belong to
a tesseract, form an orbit under the cyclic action and were already
mentioned in \cite[page~19]{Su2}. The corresponding partitions are
$(3^3)$, $(5\,2^3)$, and $(4^2\,1^3)$. The three vertices
belong to three cubes with the fixed-point vertex $g_{(4\,3\,2\,1)}$
as the opposite vertex in each of these three cubes. The following
diagram shows the three relevant subposets in $\mathbb Y_9^\circ$.
\begin{center}
\xymatrix@R=10mm@C=7mm{&&&(4\,3^2\,1)\ar@{-}[dl]\ar@{-}[d]\ar@{-}[dr]\\
&&(3^3\,1)\ar@{-}[d]\ar@{-}[dr]&\underline{(4\,3\,2\,1)}\ar@{-}[dl]\ar@{-}[dr]
&(4\,3^2)\ar@{-}[dl]\ar@{-}[d]\\
&&(3^2\,2\,1)\ar@{-}[dr]&\fbox{$(3^3)$}\ar@{-}[d]
\ar@{}[dddll]|{\mbox{\Large$\raisebox{2mm}{$\textstyle\tau$}\!\!\swarrow$}}
&(4\,3\,2)\ar@{-}[dl]\\
&&&(3^2\,2)\\
&(5\,3\,2^2)\ar@{-}[dl]\ar@{-}[d]\ar@{-}[dr]
&&&&(4^2\,2\,1^2)\ar@{-}[dl]\ar@{-}[d]\ar@{-}[dr]\\
(4\,3\,2^2)\ar@{-}[d]\ar@{-}[dr]&\fbox{$(5\,2^3)$}\ar@{-}[dl]\ar@{-}[dr]
\ar@{}[rrrr]_{
\mbox{\Large$\substack{\textstyle\longrightarrow\\[1mm]\textstyle\tau}$}}
&(5\,3\,2\,1)\ar@{-}[dl]\ar@{-}[d]
&&(4\,3\,2\,1^2)\ar@{-}[d]\ar@{-}[dr]
&\fbox{$(4^2\,1^3)$}\ar@{-}[dl]\ar@{-}[dr]
\ar@{}[uuull]|{
\mbox{\Large$\nwarrow\!\!\raisebox{2mm}{$\textstyle\tau$}$}}
&(4^2\,2\,1)\ar@{-}[dl]\ar@{-}[d]\\
(4\,2^3)\ar@{-}[dr]&\underline{(4\,3\,2\,1)}\ar@{-}[d]&(5\,2^2\,1)\ar@{-}[dl]
&&(4\,3\,1^3)\ar@{-}[dr]&\underline{(4\,3\,2\,1)}\ar@{-}[d]
&(4^2\,1^2)\ar@{-}[dl]\\
&(4\,2^2\,1)&&&&(4\,3\,1^2)
}
\end{center}
The three additional vertices and their incident edges are drawn in blue in the
next picture, and the positions are slightly distorted to get a faithful
picture of the combinatorial structure.
\begin{center}
\includegraphics[width=9cm]{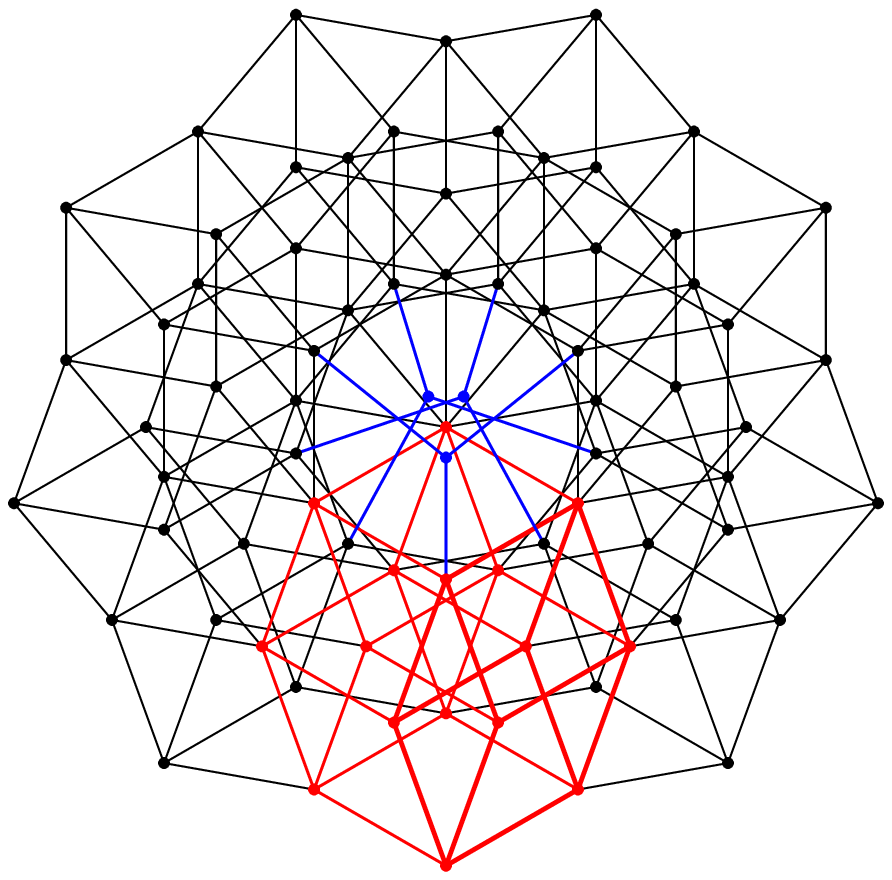}
\end{center}

The next picture displays the $1$-skeleton of $\InjHull(C_{11})$.
\begin{center}
\includegraphics[width=150mm]{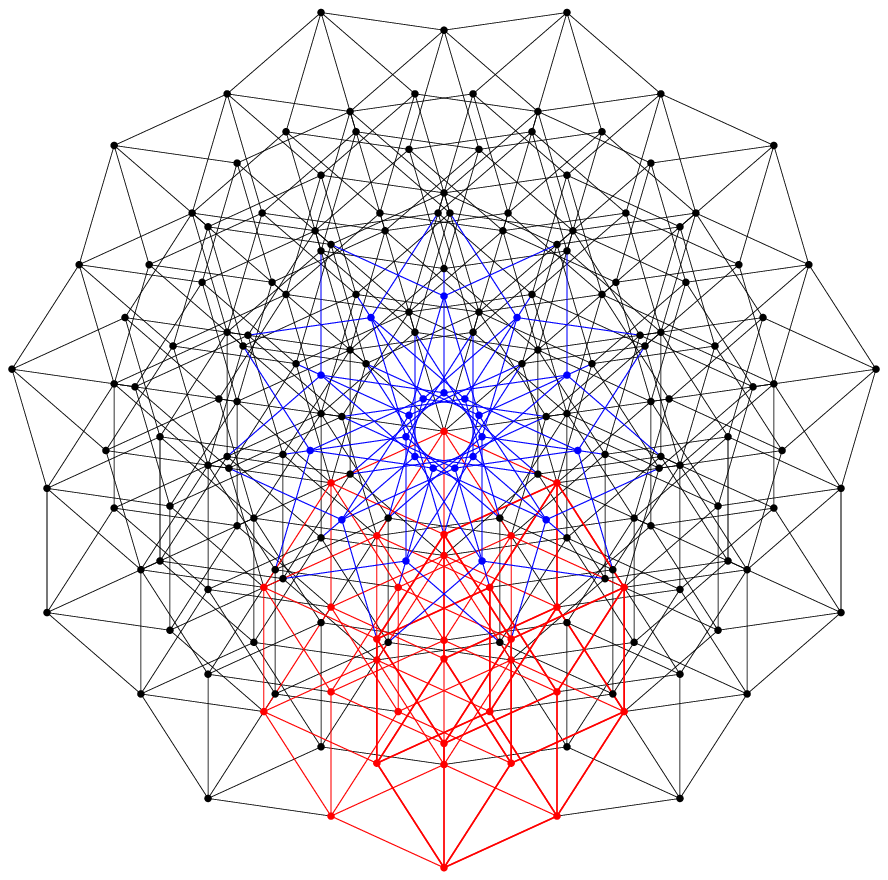}
\end{center}
There are $177$ vertices that are incident with at least one of the
eleven $5$-cubes. The $22$ additional vertices and their incident edges are
drawn in blue.

\mbox{}

It should not come as a surprise that the fibre cardinalities of
the mapping $\mathbb Y_N\to\mathbb Y_N^\circ$, $\lambda\mapsto\lambda^\circ$   
are products of Catalan numbers $C_n=\frac{1}{n+1}\binom{2n}{n}$
(a misunderstanding with the notation $C_N$ for the $N$-cycle can be excluded
from the context).
\begin{example}
For $N=11$ the fibre cardinalities are $1$, $2$, $4$, $5$, $14$, $42$, as
illustrated in the following pictures. A factor $C_n$ stands for $n$
consecutive sites in $\mathcal L_\lambda$ touching the boundary of
$\mathfrak X_N^\circ$.
\setlength{\unitlength}{0.00043333in}
\begin{center}
\begin{picture}(11124,2629)(0,-10)
\path(12,1864)(1212,664)(1512,964)
	(1812,664)(2112,964)(2412,664)
	(2712,964)(3012,664)(3312,964)
	(3612,664)(4812,1864)(4512,2164)(3312,964)
\path(4512,1564)(3912,2164)(2712,964)
\path(4212,1264)(3312,2164)(2112,964)
\path(3912,964)(2712,2164)(1512,964)
\path(3312,964)(2112,2164)(912,964)
\path(2712,964)(1512,2164)(612,1264)
\path(2112,964)(912,2164)(312,1564)
\path(1512,964)(312,2164)(12,1864)
\path(6312,1864)(7512,664)(7812,964)
	(8112,664)(8412,964)(8712,664)
	(9012,964)(9312,664)(9612,964)
	(9912,664)(11112,1864)(10812,2164)(9612,964)
\path(10812,1564)(10212,2164)(9012,964)
\path(10512,1264)(9612,2164)(8412,964)
\path(10212,964)(9012,2164)(7812,964)
\path(9612,964)(8412,2164)(7212,964)
\path(9012,964)(7812,2164)(6912,1264)
\path(8412,964)(7212,2164)(6612,1564)
\path(7812,964)(6612,2164)(6312,1864)
\texture{55888888 88555555 5522a222 a2555555 55888888 88555555 552a2a2a 2a555555 
	55888888 88555555 55a222a2 22555555 55888888 88555555 552a2a2a 2a555555 
	55888888 88555555 5522a222 a2555555 55888888 88555555 552a2a2a 2a555555 
	55888888 88555555 55a222a2 22555555 55888888 88555555 552a2a2a 2a555555 }
\shade\path(312,1564)(912,2164)(1212,1864)
	(1512,2164)(1812,1864)(2112,2164)
	(2412,1864)(2712,2164)(3012,1864)
	(3312,2164)(4212,1264)(3912,964)
	(3312,1564)(3012,1264)(2712,1564)
	(2412,1264)(2112,1564)(1812,1264)
	(1512,1564)(1212,1264)(912,1564)
	(612,1264)(312,1564)
\path(312,1564)(912,2164)(1212,1864)
	(1512,2164)(1812,1864)(2112,2164)
	(2412,1864)(2712,2164)(3012,1864)
	(3312,2164)(4212,1264)(3912,964)
	(3312,1564)(3012,1264)(2712,1564)
	(2412,1264)(2112,1564)(1812,1264)
	(1512,1564)(1212,1264)(912,1564)
	(612,1264)(312,1564)
\shade\path(6612,1564)(7212,2164)(7512,1864)
	(7812,2164)(8112,1864)(8412,2164)
	(8712,1864)(9012,2164)(9612,1564)
	(9912,1864)(10512,1264)(10212,964)
	(9912,1264)(9612,964)(9012,1564)
	(8712,1264)(8412,1564)(8112,1264)
	(7812,1564)(7512,1264)(7212,1564)
	(6912,1264)(6612,1564)
\path(6612,1564)(7212,2164)(7512,1864)
	(7812,2164)(8112,1864)(8412,2164)
	(8712,1864)(9012,2164)(9612,1564)
	(9912,1864)(10512,1264)(10212,964)
	(9912,1264)(9612,964)(9012,1564)
	(8712,1264)(8412,1564)(8112,1264)
	(7812,1564)(7512,1264)(7212,1564)
	(6912,1264)(6612,1564)
\put(312,64){\makebox(0,0)[lb]{$\lambda=(5,5,4,3,2,1)$}}
\put(6612,64){\makebox(0,0)[lb]{$\lambda=(5,4,4,3,2,1)$}}
\put(2112,2464){\makebox(0,0){$C_5=42$}}
\put(8112,2464){\makebox(0,0){$C_4=14$}}
\end{picture}
\end{center}
\begin{center}
\begin{picture}(11124,2629)(0,-10)
\path(12,1864)(1212,664)(1512,964)
	(1812,664)(2112,964)(2412,664)
	(2712,964)(3012,664)(3312,964)
	(3612,664)(4812,1864)(4512,2164)(3312,964)
\path(4512,1564)(3912,2164)(2712,964)
\path(4212,1264)(3312,2164)(2112,964)
\path(3912,964)(2712,2164)(1512,964)
\path(3312,964)(2112,2164)(912,964)
\path(2712,964)(1512,2164)(612,1264)
\path(2112,964)(912,2164)(312,1564)
\path(1512,964)(312,2164)(12,1864)
\path(6312,1864)(7512,664)(7812,964)
	(8112,664)(8412,964)(8712,664)
	(9012,964)(9312,664)(9612,964)
	(9912,664)(11112,1864)(10812,2164)(9612,964)
\path(10812,1564)(10212,2164)(9012,964)
\path(10512,1264)(9612,2164)(8412,964)
\path(10212,964)(9012,2164)(7812,964)
\path(9612,964)(8412,2164)(7212,964)
\path(9012,964)(7812,2164)(6912,1264)
\path(8412,964)(7212,2164)(6612,1564)
\path(7812,964)(6612,2164)(6312,1864)
\texture{55888888 88555555 5522a222 a2555555 55888888 88555555 552a2a2a 2a555555 
	55888888 88555555 55a222a2 22555555 55888888 88555555 552a2a2a 2a555555 
	55888888 88555555 5522a222 a2555555 55888888 88555555 552a2a2a 2a555555 
	55888888 88555555 55a222a2 22555555 55888888 88555555 552a2a2a 2a555555 }
\shade\path(312,1564)(912,2164)(1212,1864)
	(1512,2164)(1812,1864)(2112,2164)
	(2712,1564)(3012,1864)(3312,1564)
	(3612,1864)(4212,1264)(3912,964)
	(3612,1264)(3312,964)(3012,1264)
	(2712,964)(2112,1564)(1812,1264)
	(1512,1564)(1212,1264)(912,1564)
	(612,1264)(312,1564)
\path(312,1564)(912,2164)(1212,1864)
	(1512,2164)(1812,1864)(2112,2164)
	(2712,1564)(3012,1864)(3312,1564)
	(3612,1864)(4212,1264)(3912,964)
	(3612,1264)(3312,964)(3012,1264)
	(2712,964)(2112,1564)(1812,1264)
	(1512,1564)(1212,1264)(912,1564)
	(612,1264)(312,1564)
\shade\path(6612,1564)(7212,2164)(7512,1864)
	(7812,2164)(8412,1564)(9012,2164)
	(9312,1864)(9612,2164)(10512,1264)
	(10212,964)(9612,1564)(9312,1264)
	(9012,1564)(8412,964)(7812,1564)
	(7512,1264)(7212,1564)(6912,1264)(6612,1564)
\path(6612,1564)(7212,2164)(7512,1864)
	(7812,2164)(8412,1564)(9012,2164)
	(9312,1864)(9612,2164)(10512,1264)
	(10212,964)(9612,1564)(9312,1264)
	(9012,1564)(8412,964)(7812,1564)
	(7512,1264)(7212,1564)(6912,1264)(6612,1564)
\put(1512,2464){\makebox(0,0){$C_3=5$}}
\put(7512,2464){\makebox(0,0){$C_2=2$}}
\put(9312,2464){\makebox(0,0){$C_2=2$}}
\put(312,64){\makebox(0,0)[lb]{$\lambda=(5,4,3,3,2,1)$}}
\put(6612,64){\makebox(0,0)[lb]{$\lambda=(5,5,4,2,2,1)$}}
\end{picture}
\end{center}
\begin{center}
\begin{picture}(11124,2629)(0,-10)
\path(12,1864)(1212,664)(1512,964)
	(1812,664)(2112,964)(2412,664)
	(2712,964)(3012,664)(3312,964)
	(3612,664)(4812,1864)(4512,2164)(3312,964)
\path(4512,1564)(3912,2164)(2712,964)
\path(4212,1264)(3312,2164)(2112,964)
\path(3912,964)(2712,2164)(1512,964)
\path(3312,964)(2112,2164)(912,964)
\path(2712,964)(1512,2164)(612,1264)
\path(2112,964)(912,2164)(312,1564)
\path(1512,964)(312,2164)(12,1864)
\path(6312,1864)(7512,664)(7812,964)
	(8112,664)(8412,964)(8712,664)
	(9012,964)(9312,664)(9612,964)
	(9912,664)(11112,1864)(10812,2164)(9612,964)
\path(10812,1564)(10212,2164)(9012,964)
\path(10512,1264)(9612,2164)(8412,964)
\path(10212,964)(9012,2164)(7812,964)
\path(9612,964)(8412,2164)(7212,964)
\path(9012,964)(7812,2164)(6912,1264)
\path(8412,964)(7212,2164)(6612,1564)
\path(7812,964)(6612,2164)(6312,1864)
\texture{55888888 88555555 5522a222 a2555555 55888888 88555555 552a2a2a 2a555555 
	55888888 88555555 55a222a2 22555555 55888888 88555555 552a2a2a 2a555555 
	55888888 88555555 5522a222 a2555555 55888888 88555555 552a2a2a 2a555555 
	55888888 88555555 55a222a2 22555555 55888888 88555555 552a2a2a 2a555555 }
\shade\path(312,1564)(912,2164)(1212,1864)
	(1512,2164)(2112,1564)(2412,1864)
	(2712,1564)(3012,1864)(3312,1564)
	(3612,1864)(4212,1264)(3912,964)
	(3612,1264)(3312,964)(3012,1264)
	(2712,964)(2412,1264)(2112,964)
	(1512,1564)(1212,1264)(912,1564)
	(612,1264)(312,1564)
\path(312,1564)(912,2164)(1212,1864)
	(1512,2164)(2112,1564)(2412,1864)
	(2712,1564)(3012,1864)(3312,1564)
	(3612,1864)(4212,1264)(3912,964)
	(3612,1264)(3312,964)(3012,1264)
	(2712,964)(2412,1264)(2112,964)
	(1512,1564)(1212,1264)(912,1564)
	(612,1264)(312,1564)
\shade\path(6612,1564)(7212,2164)(7812,1564)
	(8112,1864)(8412,1564)(8712,1864)
	(9012,1564)(9312,1864)(9612,1564)
	(9912,1864)(10512,1264)(10212,964)
	(9912,1264)(9612,964)(9312,1264)
	(9012,964)(8712,1264)(8412,964)
	(8112,1264)(7812,964)(7212,1564)
	(6912,1264)(6612,1564)
\path(6612,1564)(7212,2164)(7812,1564)
	(8112,1864)(8412,1564)(8712,1864)
	(9012,1564)(9312,1864)(9612,1564)
	(9912,1864)(10512,1264)(10212,964)
	(9912,1264)(9612,964)(9312,1264)
	(9012,964)(8712,1264)(8412,964)
	(8112,1264)(7812,964)(7212,1564)
	(6912,1264)(6612,1564)
\put(1212,2464){\makebox(0,0){$C_2=2$}}
\put(7212,2464){\makebox(0,0){$C_1=1$}}
\put(312,64){\makebox(0,0)[lb]{$\lambda=(5,4,3,2,2,1)$}}
\put(6612,64){\makebox(0,0)[lb]{$\lambda=(5,4,3,2,1,1)$}}
\end{picture}
\end{center}
One could also take into account the exact number of the factors $C_0=1$
and $C_1=1$. The following products (meticulously listed in cyclic order)
give such a finer report for $N=11$, where each of the nineteen products
corresponds to a $\tau$-orbit in $\mathbb Y_{11}^\circ$ (in particular,
$C_0^{11}$ corresponds to the fixed point).
\begin{align*}
1&=C_0^{11}=C_1\cdot C_0^{10}=C_1\cdot C_0\cdot C_1\cdot C_0^8
=C_1\cdot C_0^2\cdot C_1\cdot C_0^7=C_1\cdot C_0^3\cdot C_1\cdot C_0^6\\
&=C_1\cdot C_0\cdot C_1\cdot C_0\cdot C_1\cdot C_0^6
=C_1\cdot C_0^2\cdot C_1\cdot C_0^3\cdot C_1\cdot C_0^3\\
2&=C_2\cdot C_0^9=C_2\cdot C_0\cdot C_1\cdot C_0^7=C_2\cdot C_0^7\cdot C_1
\cdot C_0\\
&=C_2\cdot C_0^2\cdot C_1\cdot C_0^6=C_2\cdot C_0^6\cdot C_1\cdot C_0^2
=C_2\cdot C_0^2\cdot C_1\cdot C_0^3\cdot C_1\cdot C_0^2\\
4&=C_2\cdot C_0\cdot C_2\cdot C_0^6\\
5&=C_3\cdot C_0^8=C_3\cdot C_0\cdot C_1\cdot C_0^6
=C_3\cdot C_0^6\cdot C_1\cdot C_0\\
14&=C_4\cdot C_0^7\\
42&=C_5\cdot C_0^6
\end{align*}
\end{example}

In particular, one might want to look at the subset with
fibre cardinality one
$$\mathbb Y_N^{\circ\circ}:=\bigl\{\lambda\in\mathbb Y_N^\circ\bigm|
\mathbb Y_N\ni\mu\mbox{ with }
\mu^\circ=\lambda\Longrightarrow\mu=\lambda\bigr\}.$$
The cardinality of $\mathbb Y_N^{\circ\circ}$ can be computed 
as a trace analogous to what was done for $\mathbb Y_N^\circ$.
The result for $N=2k+1$ is
\begin{align*}
\bigl|\mathbb Y_N^{\circ\circ}\bigr|&=\tr
\left(\begin{pmatrix}0&1&0\\1&1&1\\1&1&0\end{pmatrix}
\begin{pmatrix}0&1&0\\1&1&1\\1&1&0\end{pmatrix}^{k-1}
\begin{pmatrix}0&0&1\\1&1&0\\0&1&0\end{pmatrix}\right)
=[q^k]\frac{1+3q}{1-q(1+q)^2}\\
&=[q^k]\bigl(1+4 q+6 q^2+15 q^3+31 q^4+67 q^5+144 q^6+309 q^7
+664 q^8+1426 q^9+\cdots\bigr).
\end{align*}

The analogous sets for $2N$ are
$$\mathbb Y_{2N}^{\circ\circ}=
\begin{cases}
\bigl\{(N-1,N-1,N-3,N-3,\dots,2,2),\\
\phantom{\bigl\{\ }
(N,N-2,N-2,N-4,N-4,\dots,1,1)\bigr\}
&\mbox{if $N$ is odd,}\\
\varnothing&\mbox{if $N$ is even.}
\end{cases}$$
The two partitions in $\mathbb Y_{2N}^{\circ\circ}$ for $N$ odd
are interchanged by the cyclic action $\tau$. One can then look
at embeddings $\InjHull(C_{N})\subseteq\InjHull(C_{2N})$ such that
the fixed point of $\InjHull(C_{N})$ is mapped to the vertex corresponding
to either one of those $\tau^2$-fixed points in $\mathbb Y_{2N}^{\circ\circ}$.

\begin{example}
$N=5$
{\footnotesize
\begin{center}
\quad\xymatrix@R=6mm@C=-2mm{&&&&&&&&(5\,4\,3\,2\,1)\ar@{-}[dl]\ar@{-}[dr]\\
&{(4^2\,3\,2\,1)}&&{(5\,4\,3\,2)}&&&&{(5\,3^2\,2\,1)}\ar@{-}[dl]\ar@{-}[dr]
&&{(5\,4\,3\,1^2)}\ar@{-}[dl]\ar@{-}[dr]\\
{(4^2\,2^2\,1)}\ar@{-}[ur]&&{(4^2\,3\,2)}\ar@{-}[ul]\ar@{-}[ur]
&&{(5\,4\,2^2)}\ar@{-}[ul]&{\qquad\qquad}&{(4\,3^2\,2\,1)}\ar@{-}[d]
&&{(5\,3^2\,1^2)}\ar@{-}[dll]\ar@{-}[d]\ar@{-}[drr]
&&{(5\,4\,3\,1)}\ar@{-}[d]\\
{(4\,3\,2^2\,1)}\ar@{-}[u]&&{(4^2\,2^2)}\ar@{-}[ull]
\ar@{-}[u]\ar@{-}[urr]&&{(5\,4\,2\,1)}\ar@{-}[u]&&{(4\,3^2\,1^2)}\ar@{-}[dr]
&&{(5\,3\,2\,1^2)}\ar@{-}[dl]\ar@{-}[dr]&&{(5\,3^2\,1)}\ar@{-}[dl]\\
&{(4\,3\,2^2)}\ar@{-}[ul]\ar@{-}[ur]
&&{(4^2\,2\,1)}\ar@{-}[ul]\ar@{-}[ur]&&&&{(4\,3\,2\,1^2)}&&{(5\,3\,2\,1)}\\
&&{(4\,3\,2\,1)}\ar@{-}[ul]\ar@{-}[ur]}
\end{center}
}\noindent
The relevant parts of the Hasse diagram are displayed. On the left hand side
the poset structure of $\mathbb Y_5^\circ$ is realized as a subposet
in the Boolean lattice $\mathbb Y_{10}^\circ$, and on the right we have
the reversed poset structure as the image under $\tau^5$. The remaining
$2^5-2\cdot11=10$ elements of $\mathbb Y_{10}^\circ$ constitute the
$\tau$-orbit of $(5\,4\,2^2\,1)$.
\end{example}

In general, for $N=2k+1$ odd we have the embedding
$\mathbb Y_N^\circ\to\mathbb Y_{2N}^\circ$
as illustrated in the following diagram for $N=7$.
\begin{center}
\setlength{\unitlength}{0.00053333in}
\begin{picture}(10374,1539)(0,-10)
\path(312,1512)(12,1212)(1212,12)
	(1512,312)(1812,12)(2112,312)
	(2412,12)(3612,1212)(3312,1512)(2112,312)
\path(3312,912)(2712,1512)(1512,312)
\path(3012,612)(2112,1512)(912,312)
\path(2712,312)(1512,1512)(612,612)
\path(2112,312)(912,1512)(312,912)
\path(1512,312)(312,1512)
\path(4962,1062)(5862,162)(6162,462)
	(6462,162)(6762,462)(7062,162)
	(7362,462)(7662,162)(7962,462)
	(8262,162)(8562,462)(8862,162)
	(9162,462)(9462,162)(10362,1062)
	(10062,1362)(9162,462)
\path(10062,762)(9462,1362)(8562,462)
\path(9762,462)(8862,1362)(7962,462)
\path(9162,462)(8262,1362)(7362,462)
\path(8562,462)(7662,1362)(6762,462)
\path(7962,462)(7062,1362)(6162,462)
\path(7362,462)(6462,1362)(5562,462)
\path(6762,462)(5862,1362)(5262,762)
\path(6162,462)(5262,1362)(4962,1062)
\put(2712,612){\makebox(0,0){\footnotesize$\varepsilon_1$}}
\put(2112,612){\makebox(0,0){\footnotesize$\varepsilon_2$}}
\put(1512,612){\makebox(0,0){\footnotesize$\varepsilon_3$}}
\put(912,612){\makebox(0,0){\footnotesize$\varepsilon_4$}}
\put(9462,462){\makebox(0,0){\footnotesize$\varepsilon_1$}}
\put(8262,462){\makebox(0,0){\footnotesize$\varepsilon_2$}}
\put(7062,462){\makebox(0,0){\footnotesize$\varepsilon_3$}}
\put(5862,462){\makebox(0,0){\footnotesize$\varepsilon_4$}}
\put(2412,312){\makebox(0,0){\footnotesize$\delta_1$}}
\put(1812,312){\makebox(0,0){\footnotesize$\delta_2$}}
\put(1212,312){\makebox(0,0){\footnotesize$\delta_3$}}
\put(8862,462){\makebox(0,0){\footnotesize$\delta_1$}}
\put(7662,462){\makebox(0,0){\footnotesize$\delta_2$}}
\put(6462,462){\makebox(0,0){\footnotesize$\delta_3$}}
\put(4212,612){\makebox(0,0){\Large$\longmapsto$}}
\end{picture}
\end{center}
$\varepsilon_1,\delta_1,\dots,\varepsilon_k,\delta_k,\varepsilon_{k+1}\in\{0,1\}$
with the condition $\delta_i=0$ only if $\varepsilon_i=\varepsilon_{i+1}=0$
($i=1,\dots,k$), furthermore $\varepsilon_1+\varepsilon_{k+1}\leqslant1$
(this last inequality holds because otherwise the restriction regarding the
maximal hook length would be violated; but in a more uniform manner, it is
simply the condition
$\delta_{k+1}=0$ only if $\varepsilon_{k+1}=\varepsilon_{k+2}=0$ with
$\delta_{k+1}=1-\varepsilon_1$ and $\varepsilon_{k+2}=1-\delta_1$, and in fact,
look at the universal covering of the M\"obius strip to make the whole
periodicity manifest).
Let $\alpha_0=(k-1,\dots,1)\in\mathbb Y_N^\circ$ and
$\alpha_0^{(2)}=(2k,\dots,1)\in\mathbb Y_{2N}^\circ$ be the least
elements. The embedding can be written as
\begin{align*}
\mathbb Y_N^\circ\ni\lambda&=\bigl(\alpha_0\stackrel{.}{+}(
\delta_1\!+\!\varepsilon_1,\dots,\delta_{k-1}\!+\!\varepsilon_{k-1})\ 
\delta_{k}\!+\!\varepsilon_{k}\ \varepsilon_{k+1}
\bigr)\bigr|_{\textup{remove trailing zeros}}\\
&\phantom{{}=\ }\longmapsto \lambda^{(2)}:=\bigl(\alpha_0^{(2)}\stackrel{.}{+}
(\varepsilon_1,\delta_1,\dots,\varepsilon_k,\delta_k)\
\varepsilon_{k+1}\bigr)\bigr|_{\textup{remove trailing zeros}}\in\mathbb Y_{2N}^\circ.
\end{align*}
If we write $\tau_N$ and $\tau_{2N}$ for the cyclic actions, then
$\bigl(\tau_N(\lambda)\bigr)^{(2)}=(\tau_{2N})^2\bigl(\lambda^{(2)}\bigr)$.
In terms of the $\delta$- and $\varepsilon$-parameters the cyclic actions
are given by
\begin{align*}
\bigl(\delta_1,\dots,\delta_k\bigr)&\stackrel{\tau_N}{\longmapsto}
\bigl(1-\varepsilon_{k+1},\delta_1,\dots,\delta_{k-1}\bigr)\\
\bigl(\varepsilon_1,\dots,\varepsilon_{k+1}\bigr)&\stackrel{\tau_N}{\longmapsto}
\bigl(1-\delta_k,\varepsilon_1,\dots,\varepsilon_{k}\bigr)
\intertext{respectively}
\bigl(\varepsilon_1,\delta_1,\dots,\varepsilon_k,\delta_k,\varepsilon_{k+1}\bigr)
&\stackrel{\tau_{2N}}{\longmapsto}
\bigl(1-\varepsilon_{k+1},\varepsilon_1,\delta_1,\dots,
\varepsilon_k,\delta_k\bigr).
\end{align*}

One might want to look at the following continuous version
$$\mathbb Y_{2\infty+1}^\circ:=\left\{\bigl(\delta,\varepsilon\bigr):
(0,1)\times[0,1]\longrightarrow\{0,1\}\left|
\begin{array}{l}
\forall r\in(0,1):\delta(r)=0\Longrightarrow\varepsilon(r)=0\\
\varepsilon(0)=1\Longrightarrow\varepsilon(1)=0
\end{array}\right.
\right\}$$
or with $\delta,\varepsilon:\mathbb R\to\{0,1\}$
defined via $\delta(1+r)=1-\varepsilon(r)$ and $\varepsilon(1+r)=1-\delta(r)$.

Homotopically nontrivial loops of length $N$ in the M\"obius strip
$\mathfrak X_N$ parametrize the partitions in $\mathbb Y_N$, those loops
that lie in $\mathfrak X_N^\circ$ characterize the subset $\mathbb Y_N^\circ$.
More generally, one can consider those partitions in $\mathbb Y_N$
whose outer rims lie in
$$\mathfrak X_N^{(m)}:=\bigl\{(i,j)\in\mathfrak X_N\bigm|
k-m\leqslant j-i\leqslant N-k+m\bigr\}$$
where $k:=\bigl\lfloor\frac N2\bigr\rfloor$.
In particular, $\mathfrak X_N^{(1)}=\mathfrak X_N^\circ$ and
$\mathfrak X_N^{(k)}=\mathfrak X_N$. So
$$\mathbb Y_N^{(m)}:=
\bigl\{\lambda\in\mathbb Y_N\bigm|\mathcal L_\lambda\subseteq
\mathfrak X_N^{(m)}\bigr\}.$$

To count the admissible loops in the M\"obius strip, let us use
a rectangular shape instead of the triangular shape that we
employed before.

\noindent$N=2k+1$ odd$\phantom{\Bigl(}$

\setlength{\unitlength}{0.00037in}
\begin{picture}(6924,6039)(0,-10)
\path(12,5412)(612,6012)(2412,4212)
	(1812,3612)(12,5412)
\path(12,4212)(1812,6012)(2412,5412)
	(612,3612)(12,4212)
\path(1212,2412)(12,1212)(1212,12)
	(2412,1212)(1212,2412)
\path(612,612)(1812,1812)
\path(612,1812)(1812,612)
\path(4512,4812)(5712,6012)(6912,4812)
	(5712,3612)(4512,4812)
\path(5112,5412)(6312,4212)
\path(5112,4212)(6312,5412)
\path(6312,12)(6912,612)(5112,2412)
	(4512,1812)(6312,12)
\path(5112,12)(6912,1812)(6312,2412)
	(4512,612)(5112,12)
\put(312,2712){\makebox(0,0){\LARGE$\vdots$}}
\put(1212,3312){\makebox(0,0){\LARGE$\vdots$}}
\put(2112,2712){\makebox(0,0){\LARGE$\vdots$}}
\put(4812,3512){\makebox(0,0){\LARGE$\vdots$}}
\put(5712,2912){\makebox(0,0){\LARGE$\vdots$}}
\put(6612,3512){\makebox(0,0){\LARGE$\vdots$}}
\put(3312,4812){\makebox(0,0){\Large$\cdots$}}
\put(3312,1212){\makebox(0,0){\Large$\cdots$}}
\put(3312,3012){\makebox(0,0){\Large$\cdots$}}
\put(612,5412){\makebox(0,0){\footnotesize$u_0^{(0)}$}}
\put(612,4212){\makebox(0,0){\footnotesize$u_1^{(0)}$}}
\put(612,1212){\makebox(0,0){\footnotesize$u_m^{(0)}$}}
\put(1212,612){\makebox(0,0){\footnotesize$u_0^{(1)}$}}
\put(1212,1812){\makebox(0,0){\footnotesize$u_1^{(1)}$}}
\put(1212,4812){\makebox(0,0){\footnotesize$u_m^{(1)}$}}
\put(1812,5412){\makebox(0,0){\footnotesize$u_0^{(2)}$}}
\put(1812,4212){\makebox(0,0){\footnotesize$u_1^{(2)}$}}
\put(1812,1212){\makebox(0,0){\footnotesize$u_m^{(2)}$}}
\put(6312,612){\makebox(0,0){\footnotesize$u_0^{(N)}$}}
\put(6312,1812){\makebox(0,0){\footnotesize$u_1^{(N)}$}}
\put(6312,4812){\makebox(0,0){\footnotesize$u_m^{(N)}$}}
\put(5712,5412){\makebox(0,0){\footnotesize$u_0^{(2k)}$}}
\put(5712,4212){\makebox(0,0){\footnotesize$u_1^{(2k)}$}}
\put(5712,1212){\makebox(0,0){\footnotesize$u_m^{(2k)}$}}
\put(8000,3000){\makebox(0,0)[l]{$
\begin{pmatrix}
u_0^{(i+1)}\\\vdots\\\vdots\\u_m^{(i+1)}\end{pmatrix}=
\underbrace{\begin{pmatrix}
\makebox[0pt]{\raisebox{-2mm}[0pt][0pt]{\LARGE\ $0$}}
&&&1\\
&&\iddots&1\\
&\iddots&\iddots&\\
1&1&&\makebox[0pt]{\raisebox{0mm}[0pt][0pt]{\LARGE$0$\ }}
\end{pmatrix}}_{\textstyle{}=:S_m}\begin{pmatrix}
u_0^{(i)}\\\vdots\\\vdots\\u_m^{(i)}\end{pmatrix}$}}
\end{picture}\mbox{}\\[2mm]
$\bigl|\mathbb Y_{2k+1}^{(m)}\bigr|=\tr\bigl(S_m^{2k+1}\bigr)$.
In particular for $m=1$ we have the well-known expression of
the $N$th Lucas number from Remark~\ref{LucasL} as a sum of two
Fibonacci numbers:
$$L_N=|\mathbb Y_N^\circ|=\tr
\left(\begin{pmatrix}0&1\\1&1\end{pmatrix}^{\!N}\right)=
\tr\begin{pmatrix}\mathstrut F_{N-1}&F_N\\
\mathstrut  F_N&F_{N+1}\end{pmatrix}=F_{N-1}+F_{N+1}.$$

\noindent$N=2k$ even$\phantom{\Bigl(}$

\begin{picture}(6924,5139)(0,-10)
\path(12,4512)(612,5112)(2412,3312)
	(1812,2712)(12,4512)
\path(12,3312)(1812,5112)(2412,4512)
	(612,2712)(12,3312)
\path(4512,4512)(5112,5112)(6912,3312)
	(6312,2712)(4512,4512)
\path(4512,3312)(6312,5112)(6912,4512)
	(5112,2712)(4512,3312)
\path(1812,1212)(612,12)(12,612)
	(1212,1812)(2412,612)(1812,12)(612,1212)
\path(5112,1212)(6312,12)(6912,612)
	(5712,1812)(4512,612)(5112,12)(6312,1212)
\put(312,2112){\makebox(0,0){\LARGE$\vdots$}}
\put(1212,2712){\makebox(0,0){\LARGE$\vdots$}}
\put(2112,2112){\makebox(0,0){\LARGE$\vdots$}}
\put(3312,3912){\makebox(0,0){\Large$\cdots$}}
\put(3312,2112){\makebox(0,0){\Large$\cdots$}}
\put(612,4512){\makebox(0,0){\footnotesize$u_0^{(0)}$}}
\put(612,3312){\makebox(0,0){\footnotesize$u_1^{(0)}$}}
\put(4812,2112){\makebox(0,0){\LARGE$\vdots$}}
\put(5712,2712){\makebox(0,0){\LARGE$\vdots$}}
\put(6612,2112){\makebox(0,0){\LARGE$\vdots$}}
\put(3312,612){\makebox(0,0){\Large$\cdots$}}
\put(612,612){\makebox(0,0){\footnotesize$u_m^{(0)}$}}
\put(1812,4512){\makebox(0,0){\footnotesize$u_0^{(2)}$}}
\put(1812,3312){\makebox(0,0){\footnotesize$u_1^{(2)}$}}
\put(1812,612){\makebox(0,0){\footnotesize$u_m^{(2)}$}}
\put(6312,4512){\makebox(0,0){\footnotesize$u_0^{(N)}$}}
\put(6312,3312){\makebox(0,0){\footnotesize$u_1^{(N)}$}}
\put(6312,612){\makebox(0,0){\footnotesize$u_m^{(N)}$}}
\put(8000,2200){\makebox(0,0)[l]{$
\begin{pmatrix}
u_0^{(i+2)}\\[2mm]\vdots\\\vdots\\[2mm]u_m^{(i+2)}\end{pmatrix}=
\underbrace{\begin{pmatrix}
1&1&&&
\makebox[0pt]{\raisebox{-2mm}[0pt][0pt]{\LARGE$0$\ }}
\\1&2&\raisebox{0pt}[0pt][0pt]{$\ddots$}\\
&\ddots&\ddots&\ddots\\
&&\ddots&2&1\\
\makebox[0pt]{\raisebox{0mm}[0pt][0pt]{\LARGE\ $0$}}&&&1&1
\end{pmatrix}}_{\textstyle{}=:T_m}\begin{pmatrix}
u_0^{(i)}\\[2mm]\vdots\\\vdots\\[2mm]u_m^{(i)}\end{pmatrix}$}}
\end{picture}\mbox{}\\[2mm]
$\bigl|\mathbb Y_{2k}^{(m)}\bigr|=\tr\bigl(J_mT_m^k\bigr)$ where 
$J_m$ is the $(m+1)$\/$\times$\/$(m+1)$ matrix with $1$ on the antidiagonal
and $0$ elsewhere.

\subsubsection*{Acknowledgement}
It is a pleasure to thank Urs Lang for bringing to my attention
that the injective hulls of odd cycles are not well-understood.

\end{document}